\documentclass[sn-mathphys,Numbered]{sn-jnl}

\usepackage{graphicx}%
\usepackage{multirow}%
\usepackage{amsmath,amssymb,amsfonts}%
\usepackage{amsthm}%
\usepackage{mathrsfs}%
\usepackage[title]{appendix}%
\usepackage{xcolor}%
\usepackage{textcomp}%
\usepackage{manyfoot}%
\usepackage{booktabs}%
\usepackage{algorithm}%
\usepackage{algorithmicx}%
\usepackage{algpseudocode}%
\usepackage{listings}%

\theoremstyle{thmstyleone}%
\newtheorem{theorem}{Theorem}
\theoremstyle{thmstyletwo}%
\newtheorem{example}{Example}%
\newtheorem{remark}{Remark}%

\theoremstyle{thmstylethree}%
\newtheorem{definition}{Definition}%

\usepackage{mathtools}
\mathtoolsset{showonlyrefs=true}
\usepackage{bm}
\usepackage{subcaption}
\newtheorem{assumption}[theorem]{Assumption}
\newtheorem{lemma}[theorem]{Lemma}
\newtheorem{corollary}[theorem]{Corollary}
\newcommand{\e}[1]{\begin{equation}#1\end{equation}}
\newcommand{\en}[1]{\begin{equation}#1\nonumber\end{equation}}
\newcommand{\s}[1]{\begin{split}#1\end{split}}

\newcommand{\prn}[1]{\left(#1\right)}
\newcommand{\squ}[1]{\left[#1\right]}
\newcommand{\cur}[1]{\left\{#1\right\}}
\newcommand{\ang}[1]{\left\langle#1\right\rangle}
\newcommand{\nrm}[1]{\left\|#1\right\|}
\def\coloneq{\mathrel{\mathop:}=}
\def\rn{\mathbb{R}^n}
\def\rm{\mathbb{R}^m}
\def\tp{\mathrm{T}}
\def\lamx{\lambda^{A,B}_1(\boldsymbol{x})}
\def\lamxi{\lambda^{A,B}_i(\boldsymbol{x})}
\def\lamxj{\lambda^{A,B}_j(\boldsymbol{x})}
\def\lamy{\lambda^{A,B}_1(\boldsymbol{y})}
\def\cpartial{\partial}
\def\gpartial{\partial^\mathrm{GP}}
\allowdisplaybreaks

\begin{document}

\title[Article Title]{On a minimization problem of the maximum generalized eigenvalue: properties and algorithms}


\author*[1]{\fnm{Akatsuki} \sur{Nishioka}}\email{akatsuki\_nishioka@mist.i.u-tokyo.ac.jp}

\author[2]{\fnm{Mitsuru} \sur{Toyoda}}

\author[3,4]{\fnm{Mirai} \sur{Tanaka}}

\author[1,5]{\fnm{Yoshihiro} \sur{Kanno}}

\affil*[1]{\orgdiv{Department of Mathematical Informatics}, \orgname{The University of Tokyo}, \orgaddress{\street{Hongo 7‑3‑1}, \city{Bunkyo‑ku}, \postcode{113-8656}, \state{Tokyo}, \country{Japan}}}

\affil[2]{\orgdiv{Department of Mechanical Systems Engineering}, \orgname{Tokyo Metropolitan University}, \orgaddress{\street{6-6 Asahigaoka}, \city{Hino-shi}, \postcode{191-0065}, \state{Tokyo}, \country{Japan}}}

\affil[3]{\orgdiv{Department of Statistical Inference and Mathematics}, \orgname{The Institute of Statistical Mathematics}, \orgaddress{\street{10-3 Midori-cho}, \city{Tachikawa-shi}, \postcode{190-8562}, \state{Tokyo}, \country{Japan}}}

\affil[4]{\orgdiv{Continuous Optimization Team}, \orgname{RIKEN
Center for Advanced Intelligence Project}, \orgaddress{\street{Nihonbashi 1-chome Mitsui Building, 15th floor, 1-4-1}, \city{Nihonbashi, Chuo-ku}, \postcode{103-0027}, \state{Tokyo}, \country{Japan}}}

\affil[5]{\orgdiv{Mathematics and Informatics Center}, \orgname{The University of Tokyo}, \orgaddress{\street{Hongo 7‑3‑1}, \city{Bunkyo‑ku}, \postcode{113-8656}, \state{Tokyo}, \country{Japan}}}


\abstract{
We study properties and algorithms of a minimization problem of the maximum generalized eigenvalue of symmetric-matrix-valued affine functions, which is nonsmooth and quasiconvex, and has application to eigenfrequency optimization of truss structures. We derive an explicit formula of the Clarke subdifferential of the maximum generalized eigenvalue and prove the maximum generalized eigenvalue is a pseudoconvex function, which is a subclass of a quasiconvex function, under suitable assumptions. Then, we consider smoothing methods to solve the problem. We introduce a smooth approximation of the maximum generalized eigenvalue and prove the convergence rate of the smoothing projected gradient method to a global optimal solution in the considered problem. Also, some heuristic techniques to reduce the computational costs, acceleration and inexact smoothing, are proposed and evaluated by numerical experiments.
}




\maketitle

\section{Introduction}

In this paper, we consider a minimization problem of the maximum generalized eigenvalue over a nonempty compact convex set $S\subset\mathbb{R}^m_{>0}$:
\e{
\underset{\bm{x}\in S}{\mathrm{Minimize}}\ \ \lamx\coloneq\lambda_1 (A(\bm{x}),B(\bm{x})).
\label{p}
}
The function $\lamx$ returns the largest one of the generalized eigenvalues $\lambda_i\in\mathbb{R}\ (i=1,\ldots,n)$ of a pair of matrices $(A(\bm{x}),B(\bm{x}))$, which satisfy
\e{
A(\bm{x})\bm{v}_i=\lambda_i B(\bm{x})\bm{v}_i\ \ (i=1,\ldots,n),
}
where $\bm{v}_i\neq\bm{0}$ are generalized eigenvectors and $A,B:\mathbb{R}^m_{>0}\to\mathbb{S}^n$ are symmetric-matrix-valued functions satisfying $B(\bm{x})\succ0$ for any $\bm{x}\in\mathbb{R}^m_{>0}$. We assume that $A,B:\mathbb{R}^m_{>0}\to\mathbb{S}^n$ are affine functions with respect to $\bm{x}$ in our problem \eqref{p}. Complete definitions and assumptions are introduced in the following sections.

\subsection{Motivating example}

Optimization problems involving generalized eigenvalues are an important class of problems in structural optimization representing various phenomena such as vibration and buckling \cite{achtziger07smo,achtziger07siam,deaton14,ferrari19,ohsaki99,seyranian94,torii17}. The squares of eigenfrequencies of a structure (machine, building, etc.) are formulated as generalized eigenvalues of the stiffness matrix $K(\bm{x})$ and the mass matrix $M(\bm{x})+M_0$, where $\bm{x}$ denotes a variable determining the design of the structure. By maximizing the minimum eigenfrequency, we can obtain a vibration-resistant structure \cite{achtziger07smo,achtziger07siam,ohsaki99}. When we consider a truss structure, the stiffness matrix and the mass matrix become affine functions with respect to a design variable $\bm{x}$ (a vector of member cross-sectional areas), and thus we can formulate the minimum eigenfrequency maximization problem of a truss structure as problem \eqref{p} by setting $A(\bm{x})=-K(\bm{x})$ and $B(\bm{x})=M(\bm{x})+M_0$. A theoretical study of eigenfrequency optimization problems of truss structures is conducted in \cite{achtziger07siam} (see also an application-oriented paper by the same authors \cite{achtziger07smo}). Achtziger and Ko\v{c}vara \cite{achtziger07siam} showed that the maximum generalized eigenvalue is a quasiconvex function when matrix-valued functions $A,B$ are affine. They also proposed solution methods based on nonlinear semidefinite programming and the bisection method. In contrast, when we consider buckling load maximization or structural optimization of continua, matrices $A(\bm{x}),B(\bm{x})$ can be nonlinear \cite{ferrari19,kocvara02,torii17}. Such applications are not considered in this paper. In addition to structural optimization, there is another application of generalized eigenvalue optimization (with matrix variables) to control theory \cite{boyd93}.

\subsection{Related work}

There are studies of interior-point methods for the maximum generalized eigenvalue minimization problem \cite{boyd93,nesterov95}. They regard the problem as a generalized fractional programming problem since the maximum generalized eigenvalue can be written using the generalized Rayleigh quotient. Note that many optimization algorithms \cite{bot17,crouzeix91} for fractional programming are not directly applicable to our problem since its objective function is the maximum of infinitely many fraction-type functions.

Optimization problems involving standard (i.e., not generalized) eigenvalues have been studied extensively. The maximum standard eigenvalue of a symmetric-matrix-valued affine function is a nonsmooth convex function, and a minimization problem of it can be solved by linear semidefinite programming \cite{lewis96}. There are several other algorithms including the spectral bundle method \cite{helmberg00} and the smoothing method \cite{nesterov07}. These two are considered to be efficient in large-scale problems. When the symmetric-matrix-valued function is nonlinear, the maximum eigenvalue can be nonconvex. The applications include control theory \cite{lv15} and structural optimization \cite{takezawa11}. In this case, there are solution methods by the nonlinear semidefinite programming \cite{holmberg15}, the bundle method \cite{apkarian08} and the smoothing method \cite{nishioka23smao}. Generalized eigenvalues can be transformed into standard eigenvalues under some assumptions; however, the corresponding symmetric-matrix-valued function becomes nonlinear in general.

Quasiconvex optimization has applications to economics \cite{hu15} and machine learning \cite{hazan15} in addition to generalized eigenvalue optimization in structural engineering. Although quasiconvex optimization problems can have non-optimal stationary points and discontinuous points, global optimization is possible in some cases. Greenberg and Pierskalla \cite{greenberg73} introduced a specialized subdifferential (called the Greenberg--Pierskalla subdifferential in \cite{hu15}) which gives a necessary and sufficient condition for the global optimality of a quasiconvex optimization problem. Note that this subdifferential is not necessarily computable for a quasiconvex function. Solution methods for quasiconvex optimization problems include subgradient methods using variants of the Greenberg--Pierskalla subdifferential \cite{kiwiel01,hu15,yang22} and the bisection method with a convex feasibility subproblem \cite{boyd04}. Pseudoconvex functions are a subclass of quasiconvex functions. A pseudoconvex optimization problem has a nice property similar to a convex optimization problem; every Clarke stationary point of a pseudoconvex optimization problem is a global optimal solution. For properties and algorithms of pseudoconvex optimization, see \cite{bian18,soleimani07}.

For unconstrained nonsmooth convex optimization problems, it is known that the optimal iteration complexity of algorithms with the first-order oracle (objective values and subgradients) is $O(\epsilon^{-2})$ in terms of the error of the objective value, which is achieved by the subgradient method \cite{nesterov18}. This can be improved using stronger information than the first-order oracle such as a smooth approximation of the objective function. Nesterov \cite{nesterov05} proposed a smoothing method combined with an accelerated or fast gradient method with $O(\epsilon^{-1})$ complexity, which is faster than the subgradient method. It solves a smoothed problem with the fixed smoothing parameter depending on the desired accuracy $\epsilon$. Therefore, there is no asymptotic convergence guarantee to the solution of the original nonsmooth problem; it can only obtain an $\epsilon$-approximate solution. Bian \cite{bian21,bian20c} proposed a smoothing accelerated gradient method that updates a smoothing parameter at each iteration so that the convergence to the solution of the original nonsmooth problem with $O(k^{-1}\log k)$ rate is guaranteed. This method has worse complexity by $\log$ factor than \cite{nesterov05}. When a smooth approximation satisfies some assumptions, a smoothing accelerated gradient method with the convergence rate $O(k^{-1})$ is proposed \cite{tran17}. For smoothing methods in locally Lipshitz continuous nonconvex nonsmooth optimization, see \cite{chen12,zhang09,nishioka23smao}. These methods have convergence guarantees in the sense that every accumulation point of the specific subsequence is a Clarke stationary point; however, convergence rates are unknown.

\subsection{Contribution}

In this paper, we study mathematical properties and optimization algorithms of the maximum generalized eigenvalue minimization problem. Optimization problems involving generalized eigenvalues are an important class of problems in structural optimization. However, there are few mathematical studies of them, and heuristic methods are often used in the engineering communities. Although we only consider a simple problem setting, we aim to pave the way for theoretical and algorithmic studies of various generalized eigenvalue optimization problems.

First, we give an explicit formula of the Clarke subdifferential of the maximum generalized eigenvalue, and by using it, we show that the maximum generalized eigenvalue is a pseudoconvex function under some assumptions. To the best of our knowledge, a mathematical study (or so-called variational analysis \cite{rockafellar98}) of generalized eigenvalues based on convex and nonsmooth analysis of eigenvalues \cite{lewis96,overton92} and quasiconvex analysis \cite{greenberg73,kiwiel01} is new in the literature. 

Second, we consider first-order optimization algorithms for the maximum generalized eigenvalue minimization problem, looking ahead to the extension to large-scale problems such as optimization of 3D structures and continua, where the nonlinear semdifinite programming approach and the bisection method proposed in \cite{achtziger07siam} may no longer be effective due to high computational cost per iteration. We introduce a smooth approximation of the maximum generalized eigenvalue and prove the convergence rate of a smoothing projected gradient method based on \cite{bian20,zhang09} to a global optimal solution. In addition, we introduce heuristic techniques to reduce the computational costs: acceleration and inexact smoothing.

We evaluate the proposed method with heuristic techniques by numerical experiments. Our problem setting is different from the maximum generalized eigenvalue minimization problems in \cite{achtziger07siam} as we set the artificial positive lower bound of the optimization variable (besides, we only consider simpler constraints). However, we compare our solutions with the solutions without the artificial lower bound by numerical experiments and show that there are no significant differences.

\subsection{Paper organization and notation}

In section 2, we introduce preliminaries: definitions, problem settings, assumptions, and basic properties of the problem. In section 3, we study properties of the maximum generalized eigenvalue: the Clarke subdifferential and pseudoconvexity. In section 4, we study solution methods: smoothing methods and some heuristic techniques reducing the computational costs. In section 5, we report the results of numerical experiments of the eigenfrequency optimization of a truss structure. Finally, we present our concluding remarks in section 6. 

We use the following notation: $\mathbb{R}^m_{>0}$ is the set of $m$-dimensional real vectors with positive components. We denote the zero matrix by $0$. $A\succeq 0$ and $A\succ 0$ respectively denote that $A$ is positive semidefinite and positive definite. $\mathbb{S}^n$, $\mathbb{S}^n_{\succeq 0}$, and $\mathbb{S}^n_{\succ 0}$ are the sets of $n$-dimensional symmetric matrices, positive semidefinite matrices, and positive definite matrices, respectively. $[a_e]_{e=1}^m\in\rm$ denotes the vector whose $e$-th component is $a_e$. $\langle A,B \rangle = \mathrm{tr}(AB)$ denotes an inner product of matrices $A,B$. $\|\bm{x}\|$ denotes the Euclidean norm of a vector $\bm{x}\in\rm$.

\section{Preliminaries}

\subsection{Definitions}
\label{s_def}

First, we define generalized eigenvalues and generalized eigenvectors. See \cite{harville97,boyd93} for details.

\begin{definition}[Generalized eigenvalue]
We define generalized eigenvalues and generalized eigenvectors of a pair of $n$-dimensional symmetric matrices $(X,Y)$ with $Y$ positive definite by $\lambda_i\in\mathbb{R}$ and nonzero vectors $\bm{v}_i\in\rn$ satisfying
\e{
X\bm{v}_i=\lambda_i Y\bm{v}_i.
\label{ge}
}
\end{definition}

Generalized (and standard) eigenvalues are ordered decreasingly: $\lambda_1\ge\ldots\ge\lambda_n$ and $\bm{v}_i$ are normalized so that they satisfy $\bm{v}_i^\tp Y\bm{v}_j=\delta_{ij}$. The facts that there exist $n$ pairs of $\lambda_i\in\mathbb{R}$ and $\bm{v}_i\in\rn$ satisfying the above definition (including multiplicities), and that $\bm{v}_i$ can be normalized as $\bm{v}_i^\tp Y\bm{v}_j=\delta_{ij}$ are justified by the equivalence of \eqref{ge} with the standard eigenvalue problem of the symmetric matrix $Y^{-1/2}XY^{-1/2}$, which is explained in section 3. We write $\lambda_i(X,Y)$ to denote the generalized eigenvalues of $(X,Y)$ as they depend on $(X,Y)$. Although generalized eigenvectors $\bm{v}_i$ also depend on $(X,Y)$, they are not functions because the choice of $\bm{v}_i$ is not unique, thus we do not write arguments after $\bm{v}_i$.

Next, we define the Clarke subdifferential, which is used in the definition of a pseudoconvex function. A locally Lipschitz continuous function $f$ is almost everywhere differentiable by the Rademacher theorem \cite{rockafellar98}, and the Clarke subdifferential of $f$ can be defined. Note that there are several equivalent definitions. See \cite{clarke90} for details.

\begin{definition}[Clarke subdifferential]
Let $f:\rm\to\mathbb{R}$ be a locally Lipschitz continuous function. The Clarke subdifferential of $f$ at $\bm{x}\in\rm$ is defined by
\e{
\cpartial f(\bm{x})\coloneq\mathrm{conv}\{\underset{i\to\infty}{\lim}\nabla f(\bm{x}^i) \mid \underset{i\to\infty}{\lim} \bm{x}^i\to\bm{x},\ \bm{x}^i\in D_f \},
}
where $\mathrm{conv}X$ is the convex hull of a set $X$ and $D_f\subseteq\rm$ is the set of points where $f$ is differentiable, which is dense in $\rm$ by the Rademacher theorem. We call each element of $\cpartial f(\bm{x})$ a Clarke subgradient.
\end{definition}

The Clarke subdifferential coincides with the ordinary convex subdifferential when $f$ is convex, and $\cpartial f(\bm{x})=\{\nabla f(\bm{x})\}$ when $f$ is differentiable.

We also define a Clarke stationary point, which gives a first-order optimality condition for a locally Lipschitz continuous nonsmooth nonconvex optimization problem. See \cite{zhang09} for details.

\begin{definition}[Clarke stationary point]
For a nonsmooth optimization problem
\begin{equation}
    \underset{\bm{x}\in C}{\mathrm{Minimize}}\ \ f(\bm{x})\label{original}
\end{equation}
with a locally Lipschitz continuous (possibly nonsmooth) objective function $f:\rm\to\mathbb{R}$ and a nonempty closed convex feasible set $C\subset\rm$, a Clarke stationary point is a point $\bm{x}^*\in C$ satisfying
\begin{equation}
    \langle \bm{g},\bm{x}^*-\bm{x}\rangle\le0\ \ (\forall\bm{x}\in C)
    \label{stationary}
\end{equation}
for some $\bm{g}\in\partial f(\bm{x}^*)$.
\end{definition}

Finally, we define quasiconvex and pseudoconvex functions, which satisfy the following inclusions:
\en{
\text{convex}\subset\text{pseudoconvex}\subset\text{quasiconvex}.
}
See \cite{kiwiel01,penot97,soleimani07} for details.

\begin{definition}[Quasiconvex function\footnote{Not to be confused with a different notion of quasiconvexity (in the sense of Morrey) \cite{morrey52} often used in the calculus of variations. It is a property of the integrand of a functional.}]
Let $D\subset\rm$ be a nonempty convex set. A function $f:D\to\mathbb{R}$ is said to be quasiconvex if its sublevel set $
\{\bm{x}\in D \mid f(\bm{x})\le\alpha\}$ is convex for any $\alpha\in\mathbb{R}$.
\end{definition}

There is an equivalent definition: a function satisfying 
\e{
f(\lambda\bm{x}+(1-\lambda)\bm{y}) \le \max\{f(\bm{x}),f(\bm{y})\}
\label{quasi2}
}
for any $\bm{x},\bm{y}\in D$ and $\lambda\in[0,1]$ is said to be quasiconvex \cite{soleimani07}. 

\begin{definition}[Pseudoconvex function]
Let $D\subset\rm$ be a nonempty convex set. A locally Lipschitz continuous function $f:D\to\mathbb{R}$ is said to be pseudoconvex if
\e{
f(\bm{x})>f(\bm{y})\ \Rightarrow\ \langle \bm{g},\bm{y}-\bm{x} \rangle < 0
\label{pse}
}
for any $\bm{x},\bm{y}\in D$ and $\bm{g}\in\cpartial f(\bm{x})$.
\end{definition}

Pseudoconvexity of a convex function immediately follows from the inequality of convexity: $f(\bm{y})-f(\bm{x})\ge\langle \bm{g},\bm{y}-\bm{x} \rangle$. For a proof of quasiconvexity of a pseudoconvex function, see \cite{soleimani07}.

A minimization problem with a pseudoconvex objective function and a convex feasible set has a remarkable property similar to a convex optimization problem: every Clarke stationary point is a global optimal solution. It immediately follows from the definitions of a pseudoconvex function and a Clarke stationary point.

We give an example of a pseudoconvex function, which helps understand the pseudoconvexity of the maximum generalized eigenvalue.

\begin{example}
\label{ex1}
A typical example of a pseudoconvex (quasiconvex) function often treated in the literature \cite{bian18,hu15,kiwiel01} is $f_1:D\to\mathbb{R}$ defined as
\e{
f_1(\bm{x}) = \frac{a(\bm{x})}{b(\bm{x})},
\label{frac}
}
where $D\subset\rm$ is a nonempty convex set, $a:D\to\mathbb{R}$ is a convex function, and $b:D\to\mathbb{R}$ is an affine function satisfying $b(\bm{x})>0$ for any $\bm{x}\in D$.

Its sublevel set $\{\bm{x}\in D \mid f_1(\bm{x})\le\alpha\}=\{\bm{x}\in D \mid a(\bm{x})-\alpha b(\bm{x})\le 0\}$ is convex for any $\alpha\in\mathbb{R}$, and thus $f_1$ is a quasiconvex function. We assume differentiability of $a$ and $b$ for simplicity. Then, using the facts that $\varphi_{\bm{x}}(\bm{y})\coloneq a(\bm{y})-f_1(\bm{x})b(\bm{y})$ is convex with respect to $\bm{y}$, $\varphi_{\bm{x}}(\bm{x})=0$, and $\nabla_{\bm{y}}\varphi_{\bm{x}}(\bm{y})=\nabla_{\bm{y}} a(\bm{y})-f_1(\bm{x})\nabla_{\bm{y}} b(\bm{y})$,
the following inequalities hold:
\begin{align}
\ang{\nabla f_1(\bm{x}),\bm{y}-\bm{x}}
& = \ang{\frac{b(\bm{x})\nabla a(\bm{x})-a(\bm{x})\nabla b(\bm{x})}{b(\bm{x})^2},\bm{y}-\bm{x}}\\
& = \frac{1}{b(\bm{x})}\ang{\nabla\varphi_{\bm{x}}(\bm{x}),\bm{y}-\bm{x}}\\
& \le \frac{1}{b(\bm{x})}\prn{\varphi_{\bm{x}}(\bm{y})-\varphi_{\bm{x}}(\bm{x})}\\
& = \frac{1}{b(\bm{x})}\varphi_{\bm{x}}(\bm{y})\\
& = \frac{1}{b(\bm{x})}\prn{a(\bm{y})-f_1(\bm{x})b(\bm{y})}\\
& < 0\ \ (\forall\bm{y}\in D\ \mathrm{s.t.}\ f_1(\bm{x})>f_1(\bm{y})).
\end{align}
Therefore, $f_1$ is a pseudoconvex function.

The maximum generalized eigenvalue is the maximum of infinitely many rational functions as it can be written by using the Rayleigh quotient. Since a union of convex sets is convex, quasiconvexity is preserved under the max operation. Optimization of functions similar to \eqref{frac} is called fractional programming \cite{bot17,crouzeix91}.
\end{example}

\subsection{Problem setting}

We consider the maximum generalized eigenvalue minimization problem \eqref{p}:
\e{
\underset{\bm{x}\in S}{\mathrm{Minimize}}\ \ \lamx\coloneq\lambda_1 (A(\bm{x}),B(\bm{x}))
}
with the following assumptions.
\begin{assumption}
\ \\ \vspace{-6mm}
\begin{itemize}
    \item The feasible set $S\subset\mathbb{R}^m_{>0}$ is a nonempty compact convex set.
    \item Symmetric-matrix-valued functions $A:\mathbb{R}^m_{>0}\to\mathbb{S}^n$ and $B:\mathbb{R}^m_{>0}\to\mathbb{S}^n$ are affine functions: they can be written as $A(\bm{x})=A_0+\sum_{e=1}^m x_e A_e$ and $B(\bm{x})=B_0+\sum_{e=1}^m x_e B_e$, where $x_e\ (e=1,\ldots,m)$ is the $e$-th component of $\bm{x}\in\mathbb{R}^m_{>0}$ and $A_e,B_e\in\mathbb{S}^n\ (e=0,\ldots,m)$ are constant matrices.
    \item $B(\bm{x})\succ 0$ for any $\bm{x}\in S$.
\end{itemize}
\label{asm1}
\end{assumption}

\subsection{Application: eigenfrequency optimization of truss structures}

\begin{figure}[b]
  \centering
  \begin{tabular}{cc}
  \begin{minipage}[t]{0.5\hsize}
    \centering
    \includegraphics[width=3.5cm]{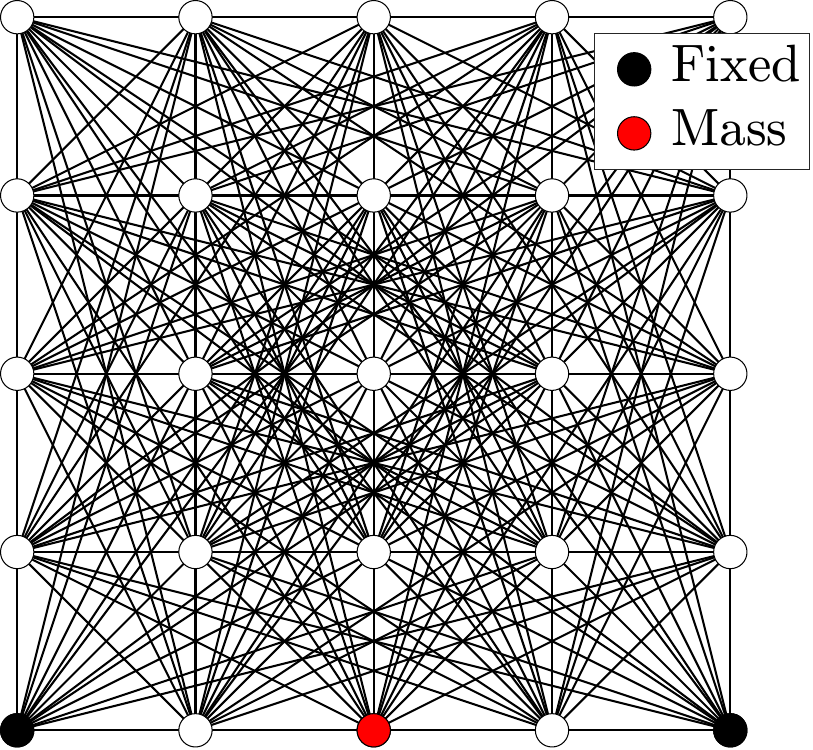}
    \vspace{-2mm}
    \subcaption{Initial design (ground structure)}
  \end{minipage} &
  \hspace{-20mm}
  \begin{minipage}[t]{0.5\hsize}
    \centering
    \includegraphics[width=3.5cm]{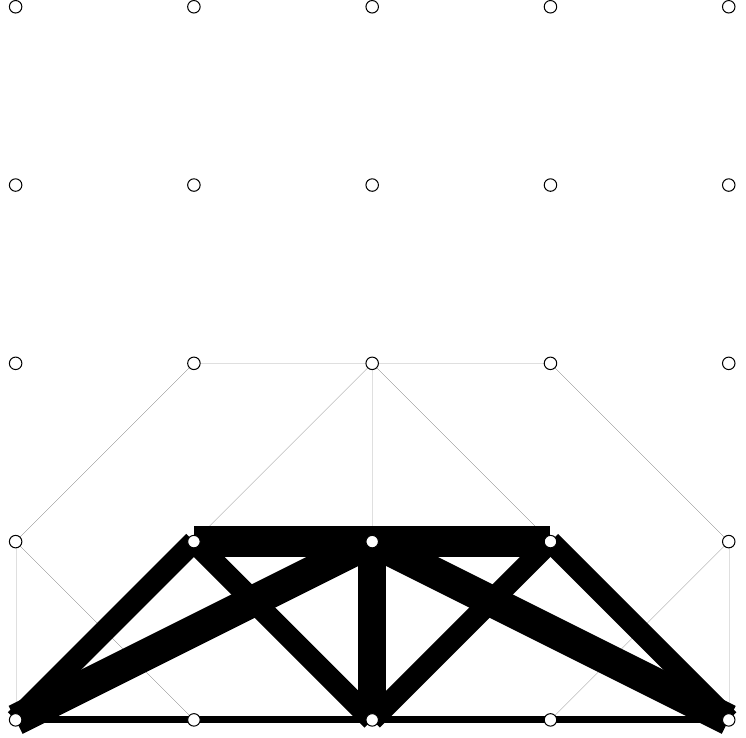}
    \vspace{-2mm}
    \subcaption{Optimal design}
  \end{minipage}
  \end{tabular}
  \caption{Eigenfrequency optimization of a truss structure}
  \label{f_setting}
\end{figure}

The minimum eigenfrequency maximization problem of a truss structure under the volume constraint \cite{achtziger07siam} shown in Fig.~\ref{f_setting} is formulated as problem \eqref{p}.
The squares of eigenfrequencies of a structure are formulated as generalized eigenvalues of the stiffness matrix $K(\bm{x})$ and the mass matrix $M(\bm{x})+M_0$. In the case of a truss structure, symmetric-matrix-valued functions $K,M:\mathbb{R}^m_{\ge 0}\to\mathbb{S}^n_{\succeq 0}$ are linear ($K_e$ and $M_e$ correspond to the element stiffness and mass matrices). A constant matrix $M_0\succeq0$ is a mass matrix of an applied non-structural mass.

By switching the sign of $K(\bm{x})$ to reformulate it as a minimization problem, the problem can be written as
\e{
\underset{\bm{x}\in S\subset\rm}{\mathrm{Minimize}}\ \ \lambda_1(-K(\bm{x}),M(\bm{x})+M_0),
\label{maxeig}
}
where
\e{
S=\{\bm{x}\in\rm \mid \bm{l}^\tp\bm{x}\le V_0,\ x_e \ge x_{\mathrm{min}}\ (e=1,\ldots,m)\}.
\label{const}
}
Here, optimization (or design) variables $x_e\ (e=1,\ldots,m)$ are cross-sectional areas of the bars of a truss structure, a constant vector $\bm{l}\in\mathbb{R}^m_{>0}$ consists of the length of each bar, $V_0>0$ is the upper bound of the volume of a truss structure, and $x_{\mathrm{min}}>0$ is the lower bound of cross-sectional areas. The first constraint in \eqref{const} is the volume constraint (or equivalently mass constraint when the density is constant), which is common in structural optimization to make an efficient structure. Therefore, the boundedness of the feasible set in Assumption \ref{asm1} is natural.

If the ground structure (see e.g., \cite{achtziger07siam}) shown in Figure \ref{f_setting}(a) is connected and the boundary is fixed, the stiffness matrix $K(\bm{x})$ is positive definite for all $\bm{x}\in[x_{\mathrm{min}},\infty)^m$ with $x_{\mathrm{min}}>0$. Moreover, since the null space of $M(\bm{x})+M_0$ is a subspace of the null space of $K(\bm{x})$ \cite{achtziger07siam}, $M(\bm{x})+M_0$ is positive definite for all $\bm{x}\in[x_{\mathrm{min}},\infty)$. Therefore, problem \eqref{maxeig} satisfies Assumption \ref{asm1}. Note that in the case that $M_0$ is positive definite (if non-structural mass is applied on every node of the ground structure), we can set $x_{\mathrm{min}}=0$ and $M(\bm{x})+M_0$ is positive definite for all $\bm{x}\in\mathbb{R}^m_{\ge 0}$.

In this paper, to ensure positive definiteness (non-singularity) of $M(\bm{x})+M_0$, we set the artificial lower bound of the optimization variables $x_{\mathrm{min}}>0$ unlike \cite{achtziger07siam}. When $x_{\mathrm{min}}=0$, the optimal structure can be topologically different from the initial structure (the ground structure), and thus the problem is called a truss topology optimization problem. Although it may not be a topology optimization problem when $x_{\mathrm{min}}>0$ to be precise, we show that the solution with a sufficiently small value such as $x_{\mathrm{min}}=10^{-8}$ can be close enough to the solution with $x_{\mathrm{min}}=0$ by numerical experiments. The extended definition and properties of the maximum eigenvalue when $x_{\mathrm{min}}=0$ (when $B(\bm{x})$ can be singular) is explained in Appendix \ref{a_dis}.

\subsection{Basic properties and problem classes}

The maximum generalized eigenvalue can be written using the generalized Rayleigh quotient \cite{boyd93} as
\e{
\lamx = \underset{\|\bm{v}\|=1}{\max}\frac{\bm{v}^\tp A(\bm{x})\bm{v}}{\bm{v}^\tp B(\bm{x})\bm{v}}.
}
It is known to be a quasiconvex function \cite{achtziger07siam}, which can be proven in a similar way to Example \ref{ex1}. Namely, a sublevel set of $\lamx$ is
\begin{align}
\{\bm{x}\in\rm \mid \lamx\le \alpha\}
& = \underset{\|\bm{v}\|=1}{\bigcap}\cur{\bm{x}\in\rm \mid \frac{\bm{v}^\tp A(\bm{x}) \bm{v}}{\bm{v}^\tp B(\bm{x}) \bm{v}} \le \alpha}\\
& = \underset{\|\bm{v}\|=1}{\bigcap}\{\bm{x}\in\rm \mid \bm{v}^\tp(A(\bm{x})-\alpha B(\bm{x}))\bm{v} \le 0\}\\
& = \{\bm{x}\in\rm \mid A(\bm{x})-\alpha B(\bm{x}) \preceq 0\},
\label{quasi}
\end{align}
which is convex. Similarly, the minimum generalized eigenvalue is a quasiconcave function ($f$ is quasiconcave iff $-f$ is quasiconvex).

Using the last equality in \eqref{quasi}, a generalized eigenvalue optimization problem can be written as a (nonlinear) semidefinite programming problem \cite{achtziger07siam}. For example, problem \eqref{p} is equivalent to
\e{
\begin{aligned}
& \underset{\bm{x}\in S,\ \lambda\in \mathbb{R}}{\mathrm{Minimize}} & & \lambda\\
& \mathrm{subject\ to} & & \lambda B(\bm{x})-A(\bm{x})\succeq 0,
\label{nsdp}
\end{aligned}
}
which is a nonlinear semidefinite programming problem because the new optimization variable $\lambda$ is multiplied by $B(\bm{x})$. In contrast, a constraint $\lamx \le c$ is convex when $c$ is constant. Therefore, the volume minimization problem of a truss structure under the minimum eigenfrequency constraint can be solved by linear semidefinite programming \cite{ohsaki99}, which is easier than \eqref{nsdp} in general.

For general quasiconvex optimization problems, we cannot construct an optimization algorithm that finds an approximate solution in a finite number of iterations; consider minimization of a quasiconvex function $f_2:\mathbb{R}\to\mathbb{R}$ defined by
\e{
f_2(x)=\begin{cases}
0 & (x=a),\\
1 & (\mathrm{otherwise}),
\end{cases}
}
where $a\in\mathbb{R}$ is a constant. Therefore, it is essential to restrict the problem class to consider convergence rates of algorithms in quasiconvex optimization. In this paper, we only consider the maximum generalized eigenvalue minimization problem \eqref{p} with Assumption \ref{asm1} and prove the convergence rate of the smoothing projected gradient method. Note that the problem \eqref{p} includes some simpler classes of problems. In fact, when $A(\bm{x})$ and $B(\bm{x})$ are diagonal matrices whose $(i,i)$-th components are affine functions $a_i(\bm{x})$ and $b_i(\bm{x})$, respectively, the problem \eqref{p} becomes a generalized linear fractional programming problem with the objective function $\max_i \{a_i(\bm{x})/b_i(\bm{x})\}$. Also, the problem \eqref{p} includes a maximum generalized eigenvalue minimization problem with matrix variables $X,Y\in\mathbb{S}^n$ and objective function $\lambda_1 (X,Y)$.

\section{Properties of the maximum generalized eigenvalue}

The maximum generalized eigenvalue is known to be a quasiconvex function under suitable assumptions \cite{achtziger07siam}. In this section, we prove that it is also a pseudoconvex function under Assumption \ref{asm1}. We first transform the maximum generalized eigenvalue into the maximum (standard) eigenvalue of a certain continuously differentiable symmetric-matrix-valued function. Then, we derive the Clarke subdifferential, which is used in the definition of pseudoconvexity, of the maximum generalized eigenvalue using the formula in \cite{overton92}. By using it, we prove the pseudoconvexity of the maximum generalized eigenvalue. Note that the positive definiteness of $B(\bm{x})$ in Assumption \ref{asm1} is essential for the proof; otherwise, the maximum generalized eigenvalue may become discontinuous \cite{achtziger07siam} and the Clarke subdifferential (hence pseudoconvexity) cannot be defined. 

\subsection{Transformation into standard eigenvalues}

In this and the next subsection, symmetric-matrix-valued functions $A,B$ are not necessarily affine. Only continuous differentiability is assumed. A generalized eigenvalue problem of a pair of matrices $(A(\bm{x}),B(\bm{x}))$ with $B(\bm{x})$ positive definite can be transformed into a standard eigenvalue problem of $B(\bm{x})^{-1/2}A(\bm{x})B(\bm{x})^{-1/2}$ by using the positive definite square root\footnote{We can set $B(\bm{x})^{1/2}=U^\tp\Lambda^{1/2}U$ since $B(\bm{x})=U^\tp\Lambda U=U^\tp\Lambda^{1/2}UU^\tp\Lambda^{1/2}U$ where $\Lambda^{1/2}$ is a diagonal matrix consisting of square roots of the eigenvalues of $B(\bm{x})$ and $U$ is a matrix consisting of the eigenvectors \cite{harville97}. Note that $B(\bm{x})^{1/2}$ is unique as shown in Lemma \ref{lmm_dif}.} $B(\bm{x})^{1/2}$ ($B(\bm{x})^{-1/2}$ denotes the inverse matrix of $B(\bm{x})^{1/2}$):
\e{
B(\bm{x})^{-1/2}A(\bm{x})B(\bm{x})^{-1/2}\bar{\bm{v}}_i=\lambda_i \bar{\bm{v}}_i\ \ (i=1,\ldots,n),
}
where $\bar{\bm{v}}_i=B(\bm{x})^{1/2}\bm{v}_i$. From this transformation, we know that $\lambda_i$ and $\bm{v}_i$ are real numbers and real vectors, respectively, that there are $n$ pairs of eigenvalues and eigenvectors including multiplicities, and that we can normalize $\bm{v}_i$ so that they satisfy $\bm{v}_i^\tp B(\bm{x})\bm{v}_i=\bar{\bm{v}}_i^\tp\bar{\bm{v}}_j=\delta_{ij}$.

We need continuous differentiability of $B(\bm{x})^{-1/2}A(\bm{x})B(\bm{x})^{-1/2}$ to apply the result on the Clarke subdifferential of the maximum standard eigenvalue shown in \cite{overton92}. See \cite{harville97} for matrix differential calculus.

\begin{lemma}
Let $A,B:\mathbb{R}^m_{>0}\to\mathbb{S}^n$ be continuously differentiable symmetric-matrix-valued functions such that $B(\bm{x})\succ 0$ for any $\bm{x}\in\mathbb{R}^m_{>0}$, and $B(\bm{x})^{1/2}$ be the positive definite square root of $B(\bm{x})$. The symmetric-matrix-valued function $B(\bm{x})^{-1/2}A(\bm{x})B(\bm{x})^{-1/2}$ is continuously differentiable on $\mathbb{R}^m_{>0}$.
\label{lmm_dif}
\end{lemma}
{\it Proof }
For continuous differentiability of the matrix product, composition, and inverse matrix, see Chapter 15 of \cite{harville97}. We prove continuous differentiability of the positive square root $\varphi:\mathbb{S}^n_{\succ 0}\to\mathbb{S}^n_{\succ 0},\ \varphi(X)=X^{1/2}$ based on the arguments in Chapter 6.1 of \cite{higham08}. We consider $\psi:\mathbb{S}^n\to\mathbb{S}^n,\ \psi(X)=X^2$ and use the inverse function theorem. Obviously, $\psi$ is continuously differentiable, and the Fr\'{e}chet derivative of $\psi$ at $X\in\mathbb{S}^n$ is, by definition, the linear mapping $D\psi(X):\mathbb{S}^n\to\mathbb{S}^n$ satisfying
\e{
\underset{\|Y\|\to 0}\lim\frac{\|\psi(X+Y)-\psi(X)-D\psi(X)Y\|}{\|Y\|}=0,
}
therefore $D\psi(X)Y=XY+YX$ for any $Y\in\mathbb{S}^n$. To apply the inverse function theorem, we show that the linear mapping $D\psi(X)$ is invertible for any $X\in\mathbb{S}^n_{\succ 0}$. Suppose that there exists a nonzero matrix $Y\in\mathbb{S}^n$ such that $D\psi(X)Y=XY+YX=0$. Then, for an eigenvalue $\lambda>0$ and the corresponding eigenvector $\bm{v}$ of $X$, we obtain $XY\bm{v}=-YX\bm{v}=-\lambda Y\bm{v}$, and thus $-\lambda<0$ is also an eigenvalue (and $Y\bm{v}$ is the corresponding eigenvector) of $X$, which contradicts to positive definiteness of $X$. Therefore, $D\psi(X)$ is invertible, and by the inverse function theorem (Theorem C.34 in \cite{lee12}), $\psi$ is a $C^1$-diffeomorphism. In other words, $\varphi$ is continuously differentiable. The uniqueness of $B(\bm{x})^{1/2}$ is also derived from the fact that $\varphi$ is an isomorphism.
\qed

\subsection{The Clarke subdifferential}

Since $B(\bm{x})^{-1/2}A(\bm{x})B(\bm{x})^{-1/2}$ is continuously differentiable and $\lamx=\lambda_1(B(\bm{x})^{-1/2}A(\bm{x})B(\bm{x})^{-1/2})$, we can apply the formula of the Clarke subdifferential of the maximum standard eigenvalue \cite{overton92} to the maximum generalized eigenvalue. Note that the maximum generalized eigenvalue $\lamx$ is locally Lipshcitz continuous because it is the composition of the Lipschitz continuous function $\lambda_1(X)$ \cite{beck17} and the continuously differentiable, hence a locally Lipschitz continuous, function $B(\bm{x})^{-1/2}A(\bm{x})B(\bm{x})^{-1/2}$.

By using the so-called spectraplex or spectrahedron (the set of matrices the eigenvalues of which belong to the unit simplex) defined by
\e{
\Upsilon_n = \{X\in\mathbb{S}^n \mid X\succeq 0,\ \mathrm{tr}(X)=1\},
}
the Clarke subdifferential of the maximum standard eigenvalue is represented as follows.

\begin{lemma}[\cite{overton92}, Theorem 3]
\label{l_ce}
Let $C:\rm\to\mathbb{S}^n$ be a continuously differentiable symmetric-matrix-valued function, $t$ be the multiplicity of the maximum (standard) eigenvalue $\lambda_1(C(\bm{x}))$, and $\bar{V}\coloneq[\bar{\bm{v}}_1,\ldots,\bar{\bm{v}}_t]\in \mathbb{R}^{n\times t}$ be a matrix consisting of the corresponding eigenvectors satisfying $\bar{V}^\tp \bar{V}=I$. Then, the Clarke subdifferential of the maximum (standard) eigenvalue $\lambda_1(C(\cdot))$ at $\bm{x}\in\rm$ can be explicitly written by
\e{
\cpartial \lambda_1(C(\bm{x})) = \left\{ \bm{g}\in\rm \mid \bm{g}=\left[\left\langle U,\bar{V}^\tp\frac{\partial C(\bm{x})}{\partial x_e}\bar{V}\right\rangle\right]_{e=1}^m,\ U\in\Upsilon_t\right\},
\label{ce}
}
where $[a_e]_{e=1}^m\in\rm$ denotes the vector whose $e$-th component is $a_e$.
\end{lemma}

We apply the above lemma to the maximum generalized eigenvalue.

\begin{theorem}[The Clarke subdifferential of the maximum generalized eigenvalue]
Let $A,B:\mathbb{R}^m_{>0}\to\mathbb{S}^n$ be continuously differentiable symmetric-matrix-valued functions such that $B(\bm{x})\succ 0$ for any $\bm{x}\in\mathbb{R}^m_{>0}$, $t$ be the multiplicity of the maximum generalized eigenvalue $\lamx$, and $V\coloneq[\bm{v}_1,\ldots,\bm{v}_t]\in \mathbb{R}^{n\times t}$ be a matrix consisting of the corresponding eigenvectors satisfying $V^\tp B(\bm{x}) V=I$. Then, the Clarke subdifferential of the maximum generalized eigenvalue $\lambda^{A,B}_1(\cdot)$ at $\bm{x}\in\mathbb{R}^m_{>0}$ can be explicitly written by
\e{
\cpartial \lamx = \left\{ \bm{g}\in\rm \mid \bm{g}=\left[\left\langle U,V^\tp\left(\frac{\partial A(\bm{x})}{\partial x_e}-\lamx \frac{\partial B(\bm{x})}{\partial x_e}\right)V\right\rangle\right]_{e=1}^m,\ U\in\Upsilon_t\right\}.
\label{eq_cla}
}
\label{t_cla}
\end{theorem}
{\it Proof }
We apply Lemma \ref{l_ce} to the maximum (standard) eigenvalue of $B(\bm{x})^{-1/2}A(\bm{x})B(\bm{x})^{-1/2}$. Let $\bar{V}\coloneq[\bar{\bm{v}}_1,\ldots,\bar{\bm{v}}_t]\in \mathbb{R}^{n\times t}$ be a matrix consisting of eigenvectors corresponding to the maximum (standard) eigenvalue of $B(\bm{x})^{-1/2}A(\bm{x})B(\bm{x})^{-1/2}$ satisfying $\bar{V}^\tp \bar{V}=I$, i.e., $\bar{V}=B(\bm{x})^{1/2}V$. We obtain
\begin{align}
& \bar{V}^\tp\left(\frac{\partial}{\partial x_e}(B(\bm{x})^{-1/2}A(\bm{x})B(\bm{x})^{-1/2})\right) \bar{V}\\
&= V^\tp B(\bm{x})^{1/2}\frac{\partial B(\bm{x})^{-1/2}}{\partial x_e}A(\bm{x})V + V^\tp\frac{\partial A(\bm{x})}{\partial x_e}V + V^\tp A(\bm{x})\frac{\partial B(\bm{x})^{-1/2}}{\partial x_e}B(\bm{x})^{1/2}V\\
&= V^\tp \frac{\partial A(\bm{x})}{\partial x_e} V - V^\tp\frac{\partial B(\bm{x})^{1/2}}{\partial x_e}B(\bm{x})^{-1/2}A(\bm{x})V - V^\tp A(\bm{x})B(\bm{x})^{-1/2}\frac{\partial B(\bm{x})^{1/2}}{\partial x_e}V\\
&= V^\tp \frac{\partial A(\bm{x})}{\partial x_e} V - \lamx V^\tp\frac{\partial B(\bm{x})^{1/2}}{\partial x_e}B(\bm{x})^{1/2}V\\
&\ \ \ - \lamx V^\tp B(\bm{x})^{1/2}\frac{\partial B(\bm{x})^{1/2}}{\partial x_e}V\\
&= V^\tp\left(\frac{\partial A(\bm{x})}{\partial x_e}-\lamx \frac{\partial B(\bm{x})}{\partial x_e}\right)V,
\label{cge2}
\end{align}
where the first equality follows by the product rule and the relation $\bar{V}=B(\bm{x})^{1/2}V$, the second equality follows by the differentiation rule of an inverse matrix $\frac{\partial F(\bm{x})^{-1}}{\partial x_e}=-F(\bm{x})^{-1}\frac{\partial F(\bm{x})}{\partial x_e}F(\bm{x})^{-1}$ \cite{harville97}, the third equality follows by the definition of generalized eigenvectors \eqref{ge} (columns of $V$ are generalized eigenvectors corresponding to the maximum generalized eigenvalue), the fourth equality follows by $\frac{\partial B(\bm{x})}{\partial x_e}=\frac{\partial}{\partial x_e}(B(\bm{x})^{1/2}B(\bm{x})^{1/2})=\frac{\partial B(\bm{x})^{1/2}}{\partial x_e}B(\bm{x})^{1/2}+B(\bm{x})^{1/2}\frac{\partial B(\bm{x})^{1/2}}{\partial x_e}$. By the relation $\lamx=\lambda_1(B(\bm{x})^{-1/2}A(\bm{x})B(\bm{x})^{-1/2})$, i.e., $\partial\lamx=\partial\lambda_1(B(\bm{x})^{-1/2}A(\bm{x})B(\bm{x})^{-1/2})$ and substitution of \eqref{cge2} into \eqref{ce} in Lemma \ref{l_ce}, we obtain \eqref{eq_cla}.
\qed

\begin{corollary}
Let $A,B:\mathbb{R}^m_{>0}\to\mathbb{S}^n$ be continuously differentiable symmetric-matrix-valued functions such that $B(\bm{x})\succ 0$ for any $\bm{x}\in\mathbb{R}^m_{>0}$. The maximum generalized eigenvalue $\lamx$ is differentiable at a point $\bm{x}\in\mathbb{R}^m_{>0}$ where its multiplicity is one, and its partial derivative is
\e{
\frac{\partial}{\partial x_e}\lamx = \bm{v}_1^\tp\left(\frac{\partial A(\bm{x})}{\partial x_e}-\lamx \frac{\partial B(\bm{x})}{\partial x_e}\right)\bm{v}_1.
}
\label{c_grad}
\end{corollary}

Note that the differentiability in Corollary \ref{c_grad} is derived from the fact that the Clarke subdifferential is a singleton \cite{clarke90}. Corollary \ref{c_grad} coincides with a result in structural optimization \cite{seyranian94} obtained by assuming the differentiability of the maximum generalized eigenvalue and differentiating both sides of the definition of generalized eigenvalues \eqref{ge}.

\subsection{Pseudoconvexity}

We prove pseudoconvexity of the maximum generalized eigenvalue under Assumption \ref{asm1} using the Clarke subdifferential obtained in Theorem \ref{t_cla}.

\begin{theorem}
Under Assumption \ref{asm1}, the maximum generalized eigenvalue $\lambda_1^{A,B}$ is pseudoconvex on $\mathbb{R}^m_{>0}$.
\label{t_pse}
\end{theorem}
{\it Proof }
Let $V\in\mathbb{R}^{m\times t}$ be the matrix defined in Theorem \ref{t_cla}. By Theorem \ref{t_cla}, a Clarke subgradient of the maximum generalized eigenvalue can be written by
\e{
\squ{\ang{U,V^\tp\left(A_e-\lamx B_e\right)V}}_{e=1}^m\in\cpartial\lamx
}
for some $U\in\Upsilon_t$. Notice that
\begin{align}
&\ang{\squ{\ang{U,V^\tp\left(A_e-\lamx B_e\right)V}}_{e=1}^m,\bm{y}-\bm{x}}\\
& = \ang{U,V^\tp\left(A(\bm{y})-\lamx B(\bm{y})\right)V}-\ang{U,V^\tp\left(A(\bm{x})-\lamx B(\bm{x})\right)V}\\
& = \ang{U,V^\tp\left(A(\bm{y})-\lamx B(\bm{y})\right)V},
\label{thr}
\end{align}
where the second equality follows from the fact that $V$ consists of generalized eigenvectors corresponding to the maximum generalized eigenvalues: $(A(\bm{x})-\lamx B(\bm{x}))V=0$. Since $\lamy<\lamx$ is equivalent to $A(\bm{y})-\lamx B(\bm{y})\prec 0$ (see \eqref{quasi}), for any $U\in\Upsilon_n$, namely for any $U=\sum_{i=1}^n \lambda_i\bm{u}_i\bm{u}_i^\tp\ (\lambda_i\ge 0,\ \sum_i\lambda_i=1,\ \|\bm{u}\|=1)$, we obtain
\begin{align}
\ang{U,V^\tp\left(A(\bm{y})-\lamx B(\bm{y})\right)V}
& = \sum_{i=1}^n\lambda_i\bm{u}_i^\tp V^\tp\left(A(\bm{y})-\lamx B(\bm{y})\right)V\bm{u}_i\\
& < 0\ \ \ \ (\forall\bm{y}\in \mathbb{R}^m_{>0}\ \mathrm{s.t.}\ \lamy<\lamx).
\end{align}
Consequently, $\lamx>\lamy \Rightarrow \langle\bm{g},\bm{y}-\bm{x}\rangle < 0$ is satisfied for any $\bm{x},\bm{y}\in\mathbb{R}^m_{>0}$ and $\bm{g}\in\partial\lamx$, and thus $\lamx$ is pseudoconvex on $\mathbb{R}^m_{>0}$.
\qed

As mentioned in Section \ref{s_def}, every Clarke stationary point of a pseudoconvex optimization problem is a global optimal solution.

\begin{corollary}
Every Clarke stationary point of the maximum generalized eigenvalue minimization problem \eqref{p} with Assumption \ref{asm1} is a global optimal solution.
\label{c_opt}
\end{corollary}

\section{Smoothing method}

In this section, we consider optimization algorithms for the maximum generalized eigenvalue minimization problem \eqref{p}. We introduce a smooth approximation of the maximum generalized eigenvalue, and then prove the convergence rate of the smoothing projected gradient method. We also introduce heuristic techniques to reduce the computational costs: acceleration and inexact smoothing.

\subsection{Smooth approximation of the maximum generalized eigenvalue}

We define a smooth approximation of the maximum generalized eigenvalue $\lamx$ by
\e{
\tilde{f}(\bm{x};\mu)\coloneq\mu\log\left(\sum_{i=1}^n \exp{\left(\frac{\lamxi}{\mu}\right)}\right),
\label{approx}
}
where $\mu>0$ is the smoothing parameter which adjusts the accuracy of the approximation and smoothness. The log-sum-exp function\footnote{It is often called the Kreisselmeier--Steinhauser (KS) function \cite{kreisselmeier80} in structural engineering.} $\mu\log\left(\sum_{i=1}^n \exp{\left(x_i/\mu\right)}\right)$ is well known as a smooth approximation of the max function $\max_i \{x_i\}$. However, as $\lambda_i^{A,B}$ is not differentiable, it is not trivial that \eqref{approx} becomes a smooth approximation of $\lambda_1^{A,B}$. In the following, the differentiability and the gradient of \eqref{approx} are derived from the theory of smooth approximations of the maximum (standard) eigenvalue \cite{beck17,nesterov07}.

\begin{lemma}
Let $A,B:\mathbb{R}^m_{>0}\to\mathbb{S}^n$ be continuously differentiable symmetric-matrix-valued functions such that $B(\bm{x})\succ 0$ for any $\bm{x}\in\mathbb{R}^m_{>0}$. For any $\mu>0$, $\tilde{f}(\cdot;\mu)$ defined by \eqref{approx} is continuously differentiable on $\mathbb{R}^m_{>0}$, and the gradient can be written as 
\e{
\nabla\tilde{f}(\bm{x};\mu)= \sum_{i=1}^n \frac{\exp(\lamxi/\mu)}{\sum_{j=1}^n \exp(\lambda^{A,B}_{j}(\bm{x})/\mu)}\left[\bm{v}_i^\tp \left(\frac{\partial A(\bm{x})}{\partial x_e}-\lamxi \frac{\partial B(\bm{x})}{\partial x_e}\right)\bm{v}_i\right]_{e=1}^m,
\label{grad}
}
where $\bm{v}_i\ (i=1,\ldots,n)$ are generalized eigenvectors corresponding to $\lamxi$, and $[a_e]_{e=1}^m\in\rm$ denotes the vector whose $e$-th component is $a_e$. In addition, for any $\bm{x}\in\mathbb{R}^m_{>0}$ and $\mu>0$, the following relation holds:
\e{
\lamx\le\tilde{f}(\bm{x};\mu)\le\lamx+\mu\log n.
\label{error}
}
\label{l_grad}
\end{lemma}
{\it Proof }
It is shown in \cite{nesterov07} that, for any $\mu>0$, the function
\e{
\mu\log\left(\sum_{i=1}^n \exp{\left(\frac{\lambda_i(X)}{\mu}\right)}\right)
\label{approx2}
}
defined on $\mathbb{S}^n$ is continuously differentiable\footnote{For differentiability of a class of functions called symmetric spectral functions, including \eqref{approx2}, see \cite{beck17,chen04,lewis01}.} for any $X\in\mathbb{S}^n$, and the gradient with respect to $X\in\mathbb{S}^n$ is 
\e{
\sum_{i=1}^n\frac{\exp(\lambda_i(X)/\mu)}{\sum_{j=1}^n \exp(\lambda_j(X)/\mu)}\bm{u}_i\bm{u}_i^\tp,
}
where $\lambda_i(X)$ and $\bm{u}_i$ denote the $i$-th eigenvalue and the corresponding normalized eigenvector of $X$, respectively. Therefore, by the relation $\lamxi =\lambda_i(B(\bm{x})^{-1/2}A(\bm{x})B(\bm{x})^{-1/2})$ and continuous differentiability of $B(\bm{x})^{-1/2}A(\bm{x})B(\bm{x})^{-1/2}$ shown in Lemma \ref{lmm_dif}, $\tilde{f}(\cdot;\mu)$ is continuously differentiable on $\mathbb{R}_{>0}^m$ for any $\mu>0$, and the partial derivative of $\tilde{f}(\cdot;\mu)$ can be written as 
\begin{align}
\frac{\partial}{\partial x_e}\tilde{f}(\bm{x};\mu)
& = \sum_{i=1}^n\frac{\exp(\lamxi/\mu)}{\sum_{j=1}^n \exp(\lamxj/\mu)}\bar{\bm{v}}_i^\tp\left(\frac{\partial}{\partial x_e}\prn{B(\bm{x})^{-1/2}A(\bm{x})B(\bm{x})^{-1/2}}\right)\bar{\bm{v}}_i\\
& = \sum_{i=1}^n\frac{\exp(\lamxi/\mu)}{\sum_{j=1}^n \exp(\lamxj/\mu)}\prn{\bm{v}_i^\tp\left(\frac{\partial A(\bm{x})}{\partial x_e}-\lamxi \frac{\partial B(\bm{x})}{\partial x_e}\right)\bm{v}_i},
\end{align}
where the last equality follows by the same argument in Equation \eqref{cge2}.

Inequality \eqref{error} immediately follows from the properties of the log-sum-exp function \cite{beck17}.
\qed

\begin{remark}
In numerical experiments, we use equivalent formulations of Equations \eqref{approx} and \eqref{grad}:
\e{
\tilde{f}(\bm{x};\mu)=\lamx+\mu\log\left(\sum_{i=1}^n \exp{\left(\frac{\lamxi-\lamx}{\mu}\right)}\right)
}
and 
\e{\s{
& \nabla\tilde{f}(\bm{x};\mu)\\
& = \sum_{i=1}^n\frac{\exp\prn{\prn{\lamxi-\lamx}/\mu}}{\sum_{j=1}^n \exp\prn{\prn{\lamxj-\lamx}/\mu}}\left[\bm{v}_i^\tp \left(\frac{\partial A(\bm{x})}{\partial x_e}-\lamxi \frac{\partial B(\bm{x})}{\partial x_e}\right)\bm{v}_i\right]_{e=1}^m
}}
to avoid overflow due to exponential terms.
\end{remark}

We show the boundedness of $\nabla\tilde{f}(\cdot;\mu)$ on a closed bounded set $S\subset\rm_{>0}$ under Assumption \ref{asm1}, which is used in the subsequent convergence analysis.

\begin{lemma}
Under Assumption \ref{asm1}, the value 
\e{
M\coloneq\max_{\bm{x}\in S}\|\nabla \tilde{f}(\bm{x};\mu)\|
\label{bound_grad}
}
is finite and independent of $\mu$, i.e., the gradient of the smooth approximation of the maximum generalized eigenvalue \eqref{grad} is bounded on a closed bounded set $S\subset\rm_{>0}$.
\end{lemma}
{\it Proof }
By the fact that $\nabla\tilde{f}(\bm{x};\mu)$ is a convex combination (Lemma \ref{l_grad}), we obtain
\e{
\|\nabla\tilde{f}(\bm{x};\mu)\|\le \max_i\nrm{\left[\bm{v}_i^\tp \left(A_e-\lamxi B_e\right)\bm{v}_i\right]_{e=1}^m}.
\label{bound}
}
Moreover, by $1=\bm{v}_i^\tp B(\bm{x})\bm{v}_i\ge\lambda_n(B(\bm{x}))\|\bm{v}_i\|^2$, we obtain $\|\bm{v}_i\|\le \max_{\bm{x}\in S}(1/\lambda_n(B(\bm{x})))$. Also, the maximum generalized eigenvalue is continuous on $S$, i.e., bounded on the bounded set $S$. Since all the terms in \eqref{bound} are bounded on $S$, $M$ is finite. Note that $\max$ in \eqref{bound_grad} is attained since $\nabla\tilde{f}(\bm{x};\mu)$ is continuous by Lemma \ref{l_grad}.
\qed

\begin{remark}
As shown in the following example, the maximum generalized eigenvalue is not Lipschitz continuous, and thus the closedness and the boundedness assumption on the feasible set $S$ is necessary to show the boundedness of the gradient of the smooth approximation.

In the example of $m=1$, $A(\bm{x})=1$, and $B(\bm{x})=x$ (i.e., $\lamx=1/x$), the gradient goes to infinity on the boundary of $\rm_{>0}$, i.e., $\{0\}$, and thus the maximum generalized eigenvalue is not Lipschitz continuous. Even if we restrict the function on $[x_\mathrm{min},\infty)^m$ with $x_\mathrm{min}>0$, we can construct an example in which the maximum generalized eigenvalue is not Lipschitz continuous. Set 
\e{
A(\bm{x})=\begin{pmatrix}
x_1 & 0 \\
0 & x_2 \\
\end{pmatrix},\ 
B(\bm{x})=\begin{pmatrix}
x_2 & 0 \\
0 & x_1 \\
\end{pmatrix},
}
then 
\begin{align}
\lamx
& = \underset{0\le v \le 1}{\max}\frac{v x_1+(1-v)x_2}{v x_2+(1-v)x_1}\\
& = \max\cur{\frac{x_1}{x_2},\frac{x_2}{x_1}},
\end{align}
where the last equality follows by the fact that the maximum is attained at an extreme point due to quasiconvexity with respect to $v$ (see \eqref{quasi2}). Let $0<x_\mathrm{min}=x_1<x_2$, then $\lamx$ is differentiable at that point and 
\e{
\frac{\partial}{\partial x_1}\lamx=-\frac{x_2}{x_\mathrm{min}^2}.
}
By taking $x_2\to\infty$, we can see the partial derivative is unbounded. Therefore, this function is not Lipschitz continuous on $[x_\mathrm{min},\infty)^m$.
\end{remark}

It is not known whether the smooth approximation \eqref{approx} is quasiconvex (pseudoconvex) nor the gradient \eqref{grad} is Lipschitz continuous.

\subsection{Convergence rate of smoothing method}

For a nonconvex nonsmooth optimization problem with a locally Lipschitz continuous objective function $f$ and a closed convex feasible set $S$, Zhang and Chen \cite{zhang09} propose the smoothing projected gradient method:
\begin{align}
& \bm{x}^{k+1}=\mathrm{\Pi}_S\left(\bm{x}^{k}-\alpha_k\nabla \tilde{f}(\bm{x}^{k};\mu_k)\right),\\
& \mu_{k+1}=
\begin{cases}
\sigma\mu_k\ \ \ \  &\text{if }\|\frac{1}{\alpha_k}(\bm{x}^{k+1}-\bm{x}^k)\|<\gamma \mu_k,\\
\mu_k &\text{otherwise},
\end{cases}
\label{spg}
\end{align}
where $\gamma>0, \sigma\in (0,1)$ are parameters, $\mathrm{\Pi}_S:\rn\to S$ is a projection operator onto $S$, and $\alpha_k>0$ is a proper stepsize (chosen by the Armijo rule with the gradient of a smooth approximation). Every accumulation point of the subsequence $\{\bm{x}^k\}_{k\in K'},\ K'=\{k\mid \mu_{k+1}=\sigma\mu_k\}$ generated by \eqref{spg} is a Clarke stationary point. Therefore, \eqref{spg} has a convergence guarantee to the global optimum in the maximum generalized eigenvalue minimization problem \eqref{p} due to pseudoconvexity. However, its convergence rate is unknown since it treats a general nonconvex nonsmooth problem. As we only consider the maximum generalized eigenvalue minimization problem \eqref{p}, which is pseudoconvex and close to a convex problem, we apply a smoothing method for a convex optimization problem to our problem \eqref{p}.

We consider the following smoothing method \cite{bian20c,bian20}:
\begin{align}
& \mu_{k}=\mu_0(k+1)^{-1/2},\\
& \alpha_{k}=\alpha_0(k+1)^{-1/2},\\
& \bm{x}^{k+1}=\mathrm{\Pi}_S\left(\bm{x}^{k}-\alpha_k\nabla \tilde{f}(\bm{x}^{k};\mu_k)\right),
\label{spgc}
\end{align}
where $\mu_0,\alpha_0>0$. The smoothing method \eqref{spgc} achieves an $O(k^{-1/2}\log k)$ convergence rate in objective values in a convex nonsmooth optimization problem \cite{bian20c,bian20}. In this paper, we prove that an $O(k^{-1/2})$ convergence rate can be achieved also in the maximum generalized eigenvalue minimization problem \eqref{p}, which is not a convex problem. Improvement of the convergence rate by $\log k$ is due to the boundedness assumption of the feasible set.

Unlike convexity, it is not straightforward to use the definition of pseudoconvexity \eqref{pse} in the convergence analysis, since it cannot connect the objective function and its gradient as an inequality. Thus, we first prove an inequality that can be used in a similar way to the inequality of convexity. Although the definition of pseudoconvexity \eqref{pse} is not used in the proof, it is based on the same approach as the proof of pseudoconvexity of the maximum generalized eigenvalue.

\begin{lemma}
For any $\bm{x},\bm{y}\in\mathbb{R}^m_{>0}$ satisfying $\lamx>\lamy$, the following inequality holds:
\e{
\langle\nabla\tilde{f}(\bm{x};\mu),\bm{y}-\bm{x}\rangle\le c_1(\bm{y})\prn{\lamy-\lamx}+c_2(\bm{y})\mu,
\label{q}
}
where $c_1(\bm{y}) \coloneq \lambda_n(B(\bm{y}))$ and $c_2(\bm{y}) \coloneq \lambda_1(B(\bm{y}))(n-1)/\mathrm{e}$.
\label{l_q}
\end{lemma}
{\it Proof }
We write $\theta_i=\exp(\lamxi/\mu)/\sum_{j=1}^n \exp(\lamxj/\mu)\ (i=1,\ldots,n)$ for simplicity. Note that $\theta_i\ge 0\ (i=1,\ldots,n)$ and $\sum_{i=1}^n\theta_i=1$. Using the formula \eqref{grad}, we obtain
\begin{align}
& \langle\nabla\tilde{f}(\bm{x};\mu),\bm{y}-\bm{x}\rangle\\
&= \sum_{i=1}^n \theta_i\left\langle\left[\bm{v}_i^\tp \left(A_e-\lamxi B_e\right)\bm{v}_i\right]_{e=1}^m,\bm{y}-\bm{x}\right\rangle\\
&= \sum_{i=1}^n \theta_i\prn{\bm{v}_i^\tp\prn{A(\bm{y})-\lamxi B(\bm{y})}\bm{v}_i-\bm{v}_i^\tp\prn{A(\bm{x})-\lamxi B(\bm{x})}\bm{v}_i}\\
&= \sum_{i=1}^n \theta_i\bm{v}_i^\tp \left(A(\bm{y})-\lamxi B(\bm{y})\right)\bm{v}_i\\
&= \sum_{i=1}^n \theta_i\bm{v}_i^\tp B(\bm{y})\bm{v}_i \left(\frac{\bm{v}_i^\tp A(\bm{y})\bm{v}_i}{\bm{v}_i^\tp B(\bm{y})\bm{v}_i}-\lamxi\right)\\
&\le \sum_{i=1}^n \theta_i\bm{v}_i^\tp B(\bm{y})\bm{v}_i \left(\lamy-\lamxi\right)\\
&= \sum_{i=1}^n \theta_i\bm{v}_i^\tp B(\bm{y})\bm{v}_i \left(\lamy-\lamx+\lamx-\lamxi\right)\\
&\le \lambda_n(B(\bm{y}))\left(\lamy-\lamx\right) + \lambda_1(B(\bm{y}))\left(\sum_{i=1}^n \theta_i\left(\lamx-\lamxi\right)\right),
\end{align}
where the last inequality follows by $\lamy-\lamx<0$ and $\lamx-\lamxi\ge 0$. Moreover,
\begin{align}
& \sum_{i=1}^n \theta_i\left(\lamx-\lamxi\right)\\
& = \sum_{i=1}^n \frac{\exp((\lamxi-\lamx/\mu)}{\sum_{j=1}^n \exp((\lamxj-\lamx/\mu)}\left(\lamx-\lamxi\right)
\end{align}
holds. Set $t_i=\lamx-\lamxi\ge 0$. By the fact that the maximum of $t\exp(-t/\mu)$ is attained at $t=\mu$, we obtain
\begin{align}
\frac{\sum_{i=1}^n t_i\exp(-t_i/\mu)}{\sum_{i=1}^n \exp(-t_i/\mu)}
& = \frac{\sum_{i=2}^n t_i\exp(-t_i/\mu)}{1+\sum_{i=2}^n \exp(-t_i/\mu)}\\
& \le \sum_{i=2}^n t_i\exp(-t_i/\mu)\\
& \le \frac{(n-1)\mu}{\mathrm{e}}.
\end{align}
Thus, the inequality \eqref{q} holds.
\qed

The inequality \eqref{q} can be interpreted as a modified version of the inequality of convexity with the coefficient $c_1(\bm{y})$ involving pseudoconvexity and the error term $c_2(\bm{y})\mu$ due to the smooth approximation. In fact, when $B(\bm{x})\equiv I$, the maximum eigenvalue and its smooth approximation become convex functions, and \eqref{q} is satisfied for any $\bm{x},\bm{y}\in\rm$ with $c_1(\bm{y})\equiv 1$. Also, when we replace $\nabla\tilde{f}(\bm{x};\mu)$ by the subgradient of the maximum generalized eigenvalue, we have $c_2(\bm{y},\mu)\equiv 0$. For the pseudoconvex maximum generalized eigenvalue, the inequality \eqref{q} is valid only for $\bm{x},\bm{y}\in\mathbb{R}^m_{>0}$ satisfying $\lamx>\lamy$. In the proof of the convergence rate, we put $\bm{y}=\bm{x}^*$, where $\bm{x}^*$ is the optimal solution.

Using Lemma \ref{l_q}, we can prove the $O\left(k^{-1/2}\right)$ convergence rate of the smoothing method \eqref{spgc} in the maximum generalized eigenvalue minimization problem \eqref{p}. The proof is similar to that of the subgradient method in convex optimization (Theorem 8.30 in \cite{beck17}).

\begin{theorem}
Let $f(\bm{x})=\lamx$ and $\{\bm{x}^k\}$ be the sequence generated by the smoothing method \eqref{spgc} for the maximum generalized eigenvalue minimization problem \eqref{p} with Assumption \ref{asm1}. Set $\alpha_{0}=1/M$, where a constant $M\coloneq\max_{\bm{x}\in S}\|\nabla \tilde{f}(\bm{x};\mu)\|$ is independent of $\mu$. The sequence $\{\bm{x}^k\}$ satisfies
\e{
\underset{j\le k}{\min}\left(f(\bm{x}^j)-f^*\right) \le \frac{2(M\Theta+(2c_2(\bm{x}^*)\mu_0+M)\log 3)}{c_1(\bm{x}^*)\sqrt{k+2}}
}
for $k\ge 2$, where $f^*$ and $\bm{x}^*$ are the global optimal value and a global optimal solution of \eqref{p}, respectively, $\Theta=\underset{\bm{x},\bm{y}\in S}{\max}\|\bm{x}-\bm{y}\|^2$, and $c_1(\cdot),c_2(\cdot)$ are defined in Lemma \ref{l_q}.
\label{t_rate}
\end{theorem}
{\it Proof }
By nonexpansiveness of projection and Lemma \ref{l_q}, we obtain
\begin{align}
\|\bm{x}^{k+1}-\bm{x}^*\|^2
& = \left\|\mathrm{\Pi}_S\left(\bm{x}^{k}-\alpha_k\nabla \tilde{f}(\bm{x}^{k};\mu_k)\right)-\bm{x}^*\right\|^2\\
& \le \left\|\bm{x}^{k}-\alpha_k\nabla \tilde{f}(\bm{x}^{k};\mu_k)-\bm{x}^*\right\|^2\\
& = \|\bm{x}^{k}-\bm{x}^*\|^2 + 2\alpha_k\langle\nabla \tilde{f}(\bm{x}^{k};\mu_k),\bm{x}^*-\bm{x}^k\rangle + \alpha_k^2\|\nabla \tilde{f}(\bm{x}^{k};\mu_k)\|^2\\
& \le \|\bm{x}^{k}-\bm{x}^*\|^2 + 2\alpha_k\prn{c_1(\bm{x}^*)(f(\bm{x}^*)-f(\bm{x}^k))+c_2(\bm{x}^*)\mu_k}\\
& \quad + \alpha_k^2\|\nabla \tilde{f}(\bm{x}^{k};\mu_k)\|^2.
\end{align}
By replacing $k$ by $i$ and summing the above inequality over $i=\lceil k/2 \rceil,\ldots,k$, we obtain
\begin{align}
& 2c_1(\bm{x}^*)\sum_{i=\lceil k/2 \rceil}^{k}\alpha_i\left(f(\bm{x}^i)-f(\bm{x}^*)\right)\\
& \le \|\bm{x}^{\lceil k/2 \rceil}-\bm{x}^*\|^2 - \|\bm{x}^{k+1}-\bm{x}^*\|^2 + \sum_{i=\lceil k/2 \rceil}^{k}\left(2c_2(\bm{x}^*)\alpha_i\mu_i + \alpha_i^2\|\nabla \tilde{f}(\bm{x}^{i};\mu_i)\|^2\right)\\
& \le \Theta + \sum_{i=\lceil k/2 \rceil}^{k}\left(\frac{2c_2(\bm{x}^*)\mu_0}{M} + 1\right)\frac{1}{i+1}.
\end{align}
Therefore, it yields
\begin{align}
\underset{j\le k}{\min}\left(f(\bm{x}^j)-f(\bm{x}^*)\right)
& \le \frac{M}{2c_1(\bm{x}^*)} \frac{\Theta+\left(\frac{2c_2(\bm{x}^*)\mu_0}{M} + 1\right)\sum_{i=\lceil k/2 \rceil}^{k}\frac{1}{i+1}}{\sum_{i=\lceil k/2 \rceil}^{k}\frac{1}{\sqrt{i+1}}}\\
& = \frac{2c_2(\bm{x}^*)\mu_0+M}{2c_1(\bm{x}^*)} \frac{\frac{M\Theta}{2c_2(\bm{x}^*)\mu_0+M}+\sum_{i=\lceil k/2 \rceil}^{k}\frac{1}{i+1}}{\sum_{i=\lceil k/2 \rceil}^{k}\frac{1}{\sqrt{i+1}}}\\
& \le \frac{2(M\Theta+(2c_2(\bm{x}^*)\mu_0+M)\log 3)}{c_1(\bm{x}^*)\sqrt{k+2}},
\end{align}
where the last inequality follows by Lemma 8.27b in \cite{beck17}:
\e{
\frac{D+\sum_{i=\lceil k/2 \rceil}^{k}\frac{1}{i+1}}{\sum_{i=\lceil k/2 \rceil}^{k}\frac{1}{\sqrt{i+1}}} \le \frac{4(D+\log 3)}{\sqrt{k+2}}\ \ (k\ge 2, D\in\mathbb{R}).
}
\qed

The above proof is based on that of the subgradient method and does not take advantage of smooth approximation. In fact, the convergence rate of the subgradient method in quasiconvex optimization \cite{kiwiel01} is the same $O(k^{-1/2})$ rate as that of Theorem \ref{t_rate}. In convex optimization, the convergence rates of the subgradient method \cite{beck17} and the smoothing method \cite{bian20,bian20c} are both $O(k^{-1/2})$.\footnote{In this argument, for simplicity, we often omit $\log k$, which appears in the rates depending on the assumptions.} In contrast, by combining the smoothing method and the accelerated gradient method, the convergence rate is improved to $O(k^{-1})$ \cite{bian20,tran17}. Therefore, we can expect to obtain the advantage of using a smooth approximation by combining the smoothing method and the accelerated gradient method also in the maximum generalized eigenvalue minimization problem \eqref{p}. 

\subsection{Heuristic techniques to reduce computational costs}

\subsubsection{Smoothing accelerated projected gradient method}

We propose the smoothing method combined with the accelerated gradient method for smooth convex optimization \cite[Algorithm 20]{daspremont21}. The proposed method lacks a convergence guarantee but converges faster than the smoothing method \eqref{spgc} and the subgradient method \cite{kiwiel01} in numerical experiments. The update formula of the smoothing accelerated projected gradient method is as follows:
\begin{align}
& \mu_{k} = \mu_0(k+1)^{-1},\\
& \alpha_{k} = \alpha_0(k+1)^{-1},\\
& \bm{y}^k = \prn{1-\frac{1}{a_k}}\bm{x}^k+\frac{1}{a_k}\bm{z}^k,\\
& \bm{z}^{k+1} = \mathrm{\Pi}_S\prn{\bm{z}^k-a_k\alpha_k\nabla \tilde{f}(\bm{y}^{k};\mu_k)},\\
& \bm{x}^{k+1} = \prn{1-\frac{1}{a_k}}\bm{x}^k+\frac{1}{a_k}\bm{z}^{k+1},\\
& a_{k+1} = \frac{1+\sqrt{4a_{k}^2+1}}{2},
\label{sapg}
\end{align}
where $\bm{x}^0=\bm{z}^0\in S$, $a_0=1$ and $\alpha_0,\mu_0>0$. All the sequences $\{\bm{x}^k\}$, $\{\bm{y}^k\}$, and $\{\bm{z}^k\}$ of \eqref{sapg} are feasible since the formulas for $\bm{y}^k$ and $\bm{x}^{k+1}$ are convex combinations. This property is important because structural optimization problems, including our problem \eqref{p}, often have objective functions only defined on the feasible set. Thus, we choose the accelerated gradient method \cite[Algorihtm 20]{daspremont21} rather than Nesterov's one \cite{nesterov83}.

The difficulty of the convergence analysis of the smoothing accelerated gradient method in our problem \eqref{p} is due to the fact that the inequality \eqref{q} has the extra coefficient $c_1(\bm{y})$ and cannot be used for any $\bm{x}$ and $\bm{y}$, unlike the inequality of convexity. There are several drawbacks of our algorithms including lack of convergence guarantee, stopping criteria, and stepsize strategy. Further studies are needed for practical applications.

\subsubsection{Inexact smoothing method}

The smooth approximation \eqref{approx} requires all the generalized eigenvalues for computation, and thus computationally costly when the size of matrices $A(\bm{x})$ and $B(\bm{x})$ is large. We propose the inexact smoothing method which replaces the gradient of the smooth approximation $\nabla\tilde{f}(\bm{x};\mu)$ in \eqref{spgc} and \eqref{sapg} by
\e{
\bm{g}_l(\bm{x};\mu)= \sum_{i=1}^l \frac{\exp(\lamxi/\mu)}{\sum_{j=1}^l \exp(\lamxj/\mu)}\left[\bm{v}_i^\tp \left(A_e-\lamxi)B_e\right)\bm{v}_i\right]_{e=1}^m
\label{inexact}
}
with $l\le n$. Namely, we only use $l$ largest generalized eigenvalues. Note that if we only use a part of generalized eigenvalues, the inexact smooth approximation $\tilde{f}_l(\bm{x};\mu)\coloneq\mu\ln\prn{\sum_{i=1}^l \exp\prn{\lamxi/\mu}}$ is may not be differentiable, and $\bm{g}_l(\bm{x};\mu)$ may not be the gradient of $\tilde{f}_l(\bm{x};\mu)$. Practically, however, the influence of generalized eigenvalues which are sufficiently smaller than the maximum one on the gradient \eqref{grad} is very small due to exponential terms. Therefore, we can expect that \eqref{inexact} is accurate enough if we use generalized eigenvalues which are sufficiently close to the maximum one. Numerical experiments show that we can take $l$ much less than $n$.

In the smoothing method without acceleration \eqref{spgc}, we can still prove the same convergence rate in Theorem \ref{t_rate} if we replace the gradient $\nabla\tilde{f}(\bm{x};\mu)$ by the inexact one $\bm{g}_l(\bm{x};\mu)$ since a similar inequality to the one in Lemma \ref{l_q} still holds. Note that the inexact smoothing method with $l=1$ corresponds to the subgradient method. We summarize these results as the following corollaries. The proofs are omitted as they are the same as the proofs of Lemma \ref{l_q} and Theorem \ref{t_rate}.

\begin{corollary}
For any $\bm{x}\in\mathbb{R}^m_{>0}$ and $\bm{y}\in\mathbb{R}^m_{>0}$ satisfying $\lamx>\lamy$, the following inequality holds:
\e{
\langle\bm{g}_l(\bm{x};\mu),\bm{y}-\bm{x}\rangle\le c_1(\bm{y})\prn{\lamy-\lamx}+c_2'(\bm{y}),
}
where $c_1(\bm{y})\coloneq\lambda_n(B(\bm{y}))$ and $c_2'(\bm{y})\coloneq\lambda_1(B(\bm{y}))(l-1)\frac{\mu}{\mathrm{e}}$.
\label{c_inexact}
\end{corollary}

\begin{corollary}
Let $f(\bm{x})=\lamx$ and $\{\bm{x}^k\}$ be the sequence generated by the smoothing method \eqref{spgc} using $\bm{g}_l(\bm{x};\mu)$ in \eqref{inexact} instead of the gradient $\nabla\tilde{f}(\bm{x};\mu)$ for the maximum generalized eigenvalue minimization problem \eqref{p} with Assumption \ref{asm1}. Set $\alpha_{0}=1/M$, where a constant $M\coloneq\max_{\bm{x}\in S}\|\nabla \tilde{f}(\bm{x};\mu)\|$ is independent of $\mu$. The sequence $\{\bm{x}^k\}$ satisfies
\e{
\underset{j\le k}{\min}\left(f(\bm{x}^j)-f^*\right) \le \frac{2(M\Theta+(2c_2'(\bm{x}^*)\mu_0+M)\log 3)}{c_1(\bm{x}^*)\sqrt{k+2}}
}
for $k\ge 2$, where $f^*$ and $\bm{x}^*$ are the global optimal value and a global optimal solution of \eqref{p}, and $\Theta=\underset{\bm{x},\bm{y}\in S}{\max}\|\bm{x}-\bm{y}\|^2$.
\label{c_inexact_rate}
\end{corollary}

For the smoothing accelerated projected gradient method \eqref{sapg} using $\bm{g}_l(\bm{x};\mu)$ instead of $\nabla\tilde{f}(\bm{x};\mu)$, the convergence is obviously not guaranteed.

\section{Numerical results}

\subsection{Setting}

We consider the eigenfrequency optimization problem of a truss structure \eqref{maxeig} and compare the following algorithms:

\begin{itemize}
\item S-APG: smoothing accelerated projected gradient method.
\item Inexact S-APG: S-APG with inexact smoothing \eqref{inexact}.
\item S-PG: smoothing projected gradient method.
\item Subgrad: subgradient method \cite{kiwiel01} (see Appendix \ref{a_sub} for details).
\end{itemize}

We set nodes, fixed points (supports), and non-structural mass as shown in Figure \ref{f_setting}. In this setting, numerical experiments show that the two maximum eigenvalues coincide near the optimal solution, and smooth optimization algorithms are not expected to work (see discussions in \cite{nishioka23smao,thore22} for example). We set the parameters of the problem as follows; the number of the optimization variables is $m=200$, the size of the matrices is $n=46$, Young's modulus of the material used in the stiffness matrix is $200$ GPa, the density of the material used in the mass matrix is $7.86\times10^{3}$ kg/m$^3$, the mass of the non-structural mass is $10^7$ kg, the distance between the nearest nodes is $1$ m, the upper limit of the volume is $V_0=0.1$ m$^3$, and the lower bound of the cross-sectional areas is $x_\mathrm{min}=10^{-8}$ m$^2$. We set the parameters of the algorithms as follows: the initial stepsizes for S-APG and Inexact S-APG are $\alpha_0=2\times10^{-6}$, the initial stepsize for S-PG is $\alpha_0=2\times10^{-7}$, the initial stepsize\footnote{A larger stepsize leads to a more oscillatory sequence in S-PG as $\alpha_k$ decays slower $O(k^{-1/2})$ rate. Also, a stepsize of Subgrad is larger since the gradient is normalized.} for Subgrad is $\alpha_0=10^{-3}$, and the initial value of the smoothing parameter is $\mu_0=10$. All the experiments have been conducted on MacBook Pro (2019, 1.4 GHz Quad-Core Intel Core i5, 8 GB memory) and MATLAB R2022b.

\subsection{Effectiveness of acceleration}

\begin{figure}[ht]
    \centering
    \includegraphics[width=8cm]{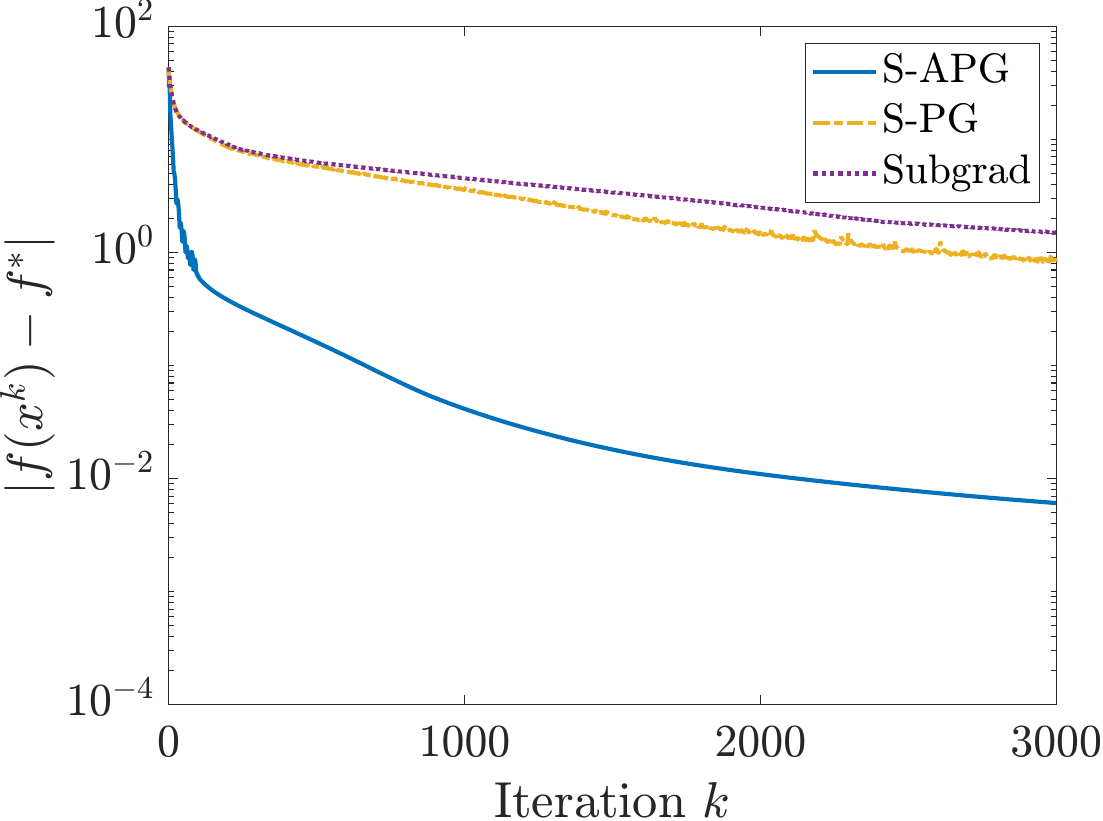}
    \caption{Difference of the objective value and the optimal value at each iteration (comparison with acceleration)}
    \label{f_obj}
\end{figure}

\begin{figure}[ht]
  \centering
  \begin{tabular}{ccc}
  \begin{minipage}[t]{0.3\hsize}
    \centering
    \includegraphics[width=3cm]{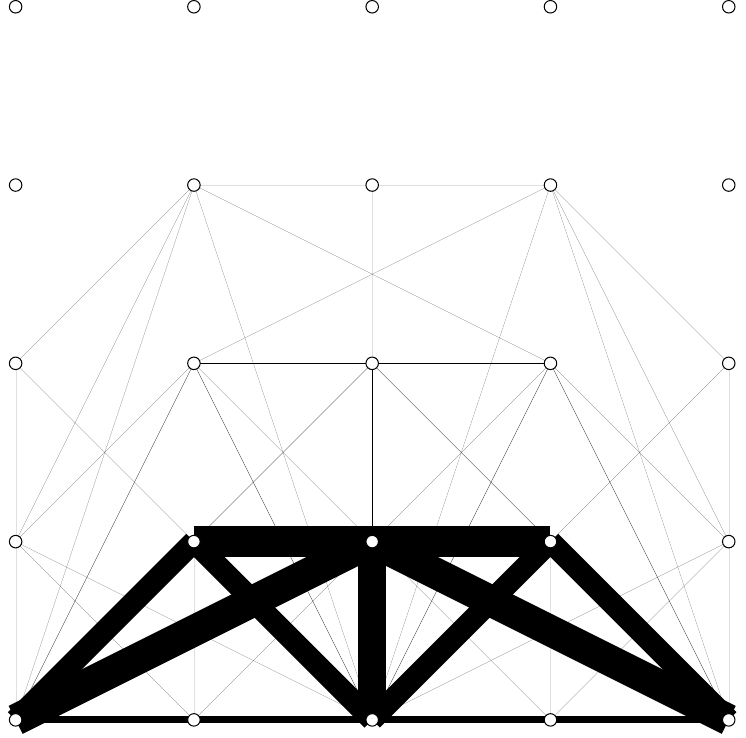}
    \vspace{-2mm}
    \subcaption{S-APG}
  \end{minipage} &
  \hspace{-4mm}
  \begin{minipage}[t]{0.3\hsize}
    \centering
    \includegraphics[width=3cm]{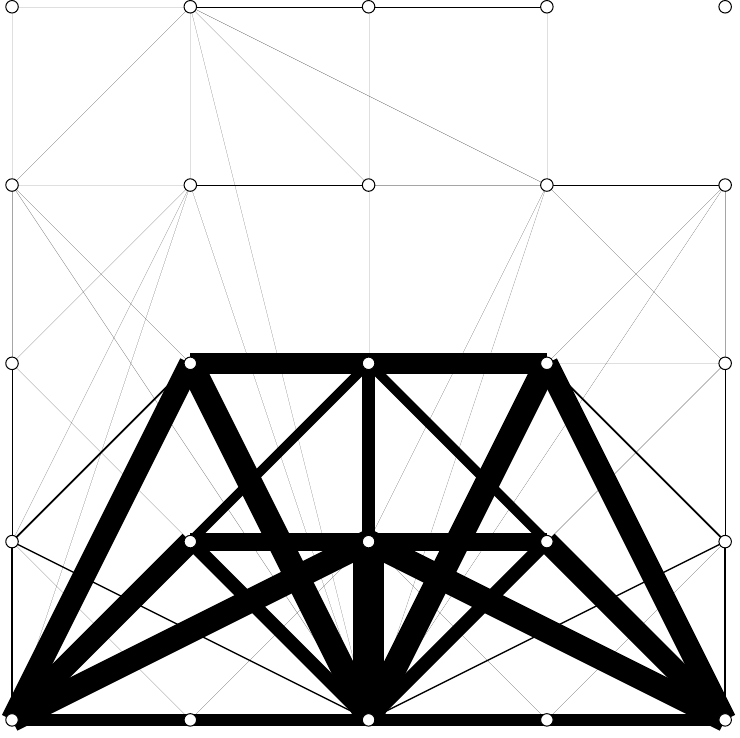}
    \vspace{-2mm}
    \subcaption{S-PG}
  \end{minipage} &
  \hspace{-4mm}
  \begin{minipage}[t]{0.3\hsize}
    \centering
    \includegraphics[width=3cm]{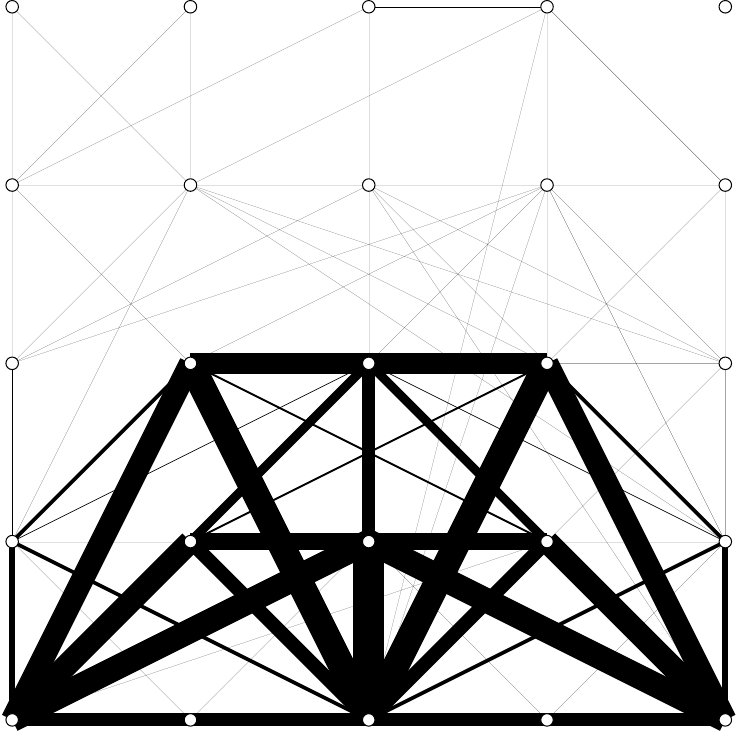}
    \vspace{-2mm}
    \subcaption{Subgrad}
  \end{minipage}
  \end{tabular}
  \caption{Designs after $3000$ iterations}
  \label{f_des}
\end{figure}

To show the effectiveness of acceleration, we compare S-APG, S-PG, and Subgrad. Figure \ref{f_obj} shows the differences between objective values $f^*$ and the optimal value in 3000 iterations, and Figure \ref{f_des} shows the truss designs after 3000 iterations. The optimal value is computed in advance by the bisection method \cite{achtziger07siam} explained in Appendix \ref{a_bis}. Bars with cross-sectional areas less than $1.5\times x_\mathrm{min}=1.5\times 10^{-8}$ are not displayed in Figure \ref{f_des}.

Figure \ref{f_obj} shows that S-PG and Subgrad converge at the same rate, which is consistent with the theoretical results. In contrast, S-APG converges faster as shown in Figures \ref{f_obj} and \ref{f_des}. We can expect that acceleration is effective in the maximum generalized eigenvalue minimization problem even though it is not a convex problem. Note that the remaining thin bars in Figure \ref{f_des} are natural because they prevent unstable nodes (as also seen in the literature \cite{achtziger07siam,ohsaki99}). Designs of S-PG and Subgrad in Figure \ref{f_des} have not converged yet.

\subsection{Effectiveness of inexact smoothing}

\begin{figure}[ht]
    \centering
    \includegraphics[width=8cm]{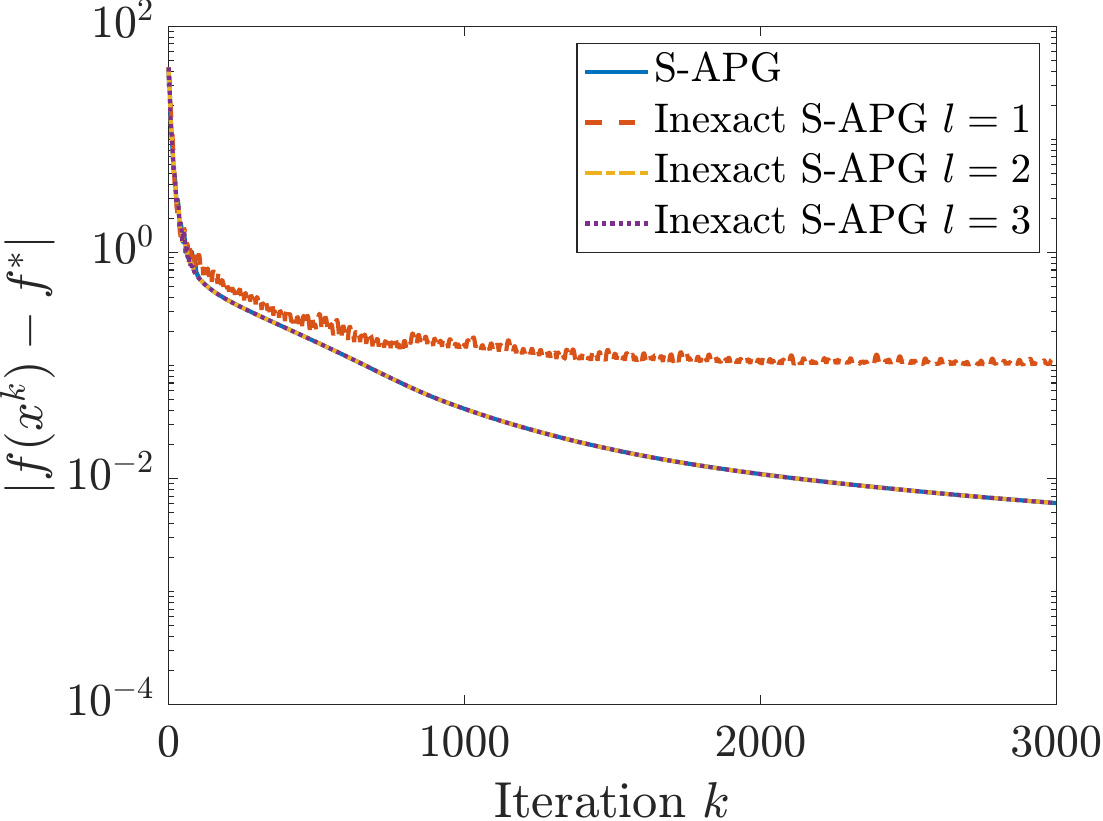}
    \caption{Difference of the objective value and the optimal value at each iteration (comparison with inexact smoothing)}
    \label{f_obj_inexact}
\end{figure}

\begin{figure}[ht]
  \centering
  \begin{tabular}{ccc}
  \begin{minipage}[t]{0.3\hsize}
    \centering
    \includegraphics[width=3cm]{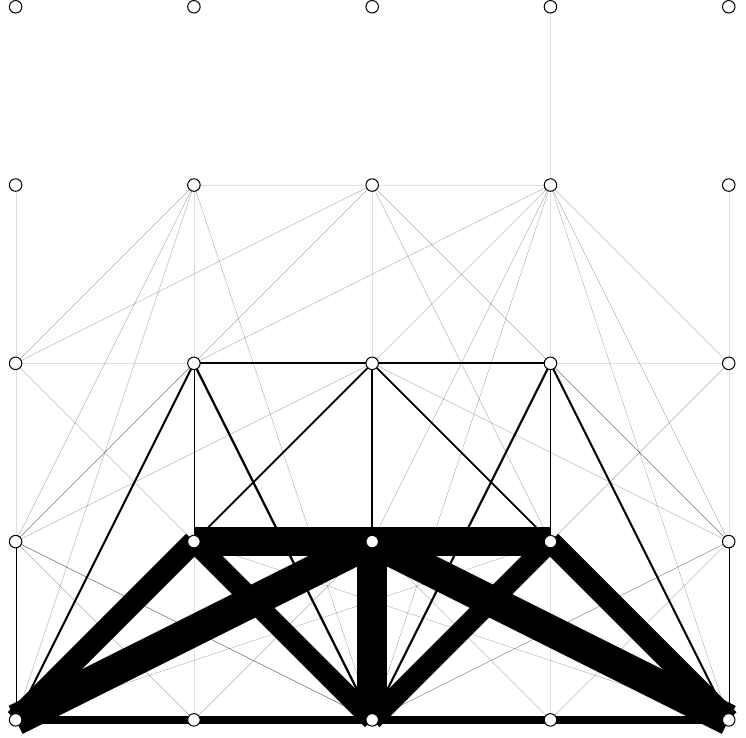}
    \vspace{-2mm}
    \subcaption{Inexact ($l=1$)}
  \end{minipage} &
  \hspace{-4mm}
  \begin{minipage}[t]{0.3\hsize}
    \centering
    \includegraphics[width=3cm]{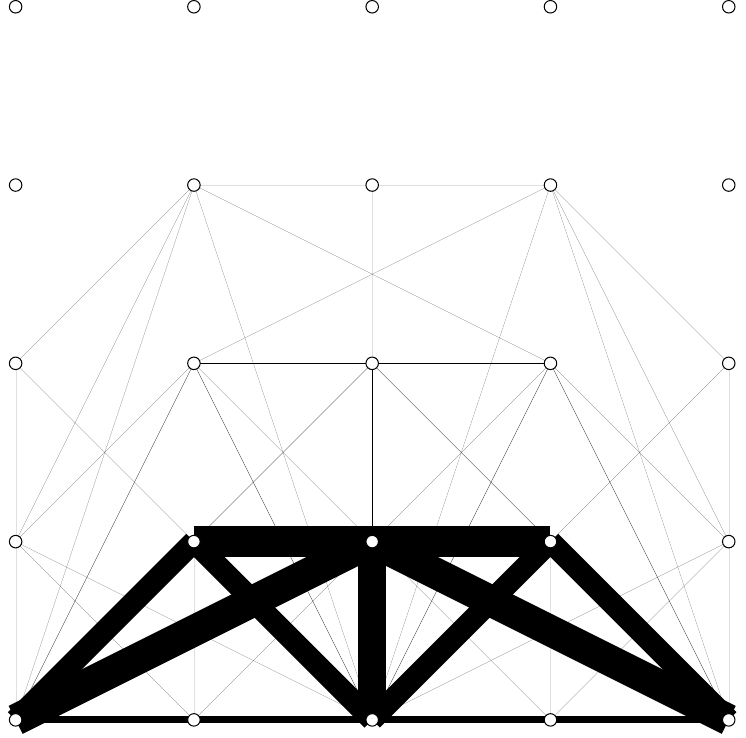}
    \vspace{-2mm}
    \subcaption{Inexact ($l=2$)}
  \end{minipage} &
  \hspace{-4mm}
  \begin{minipage}[t]{0.3\hsize}
    \centering
    \includegraphics[width=3cm]{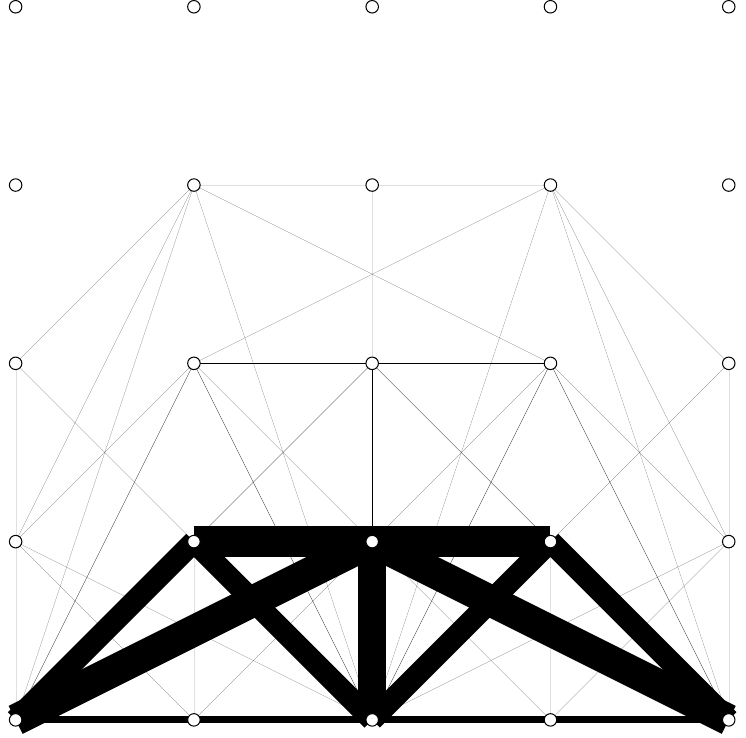}
    \vspace{-2mm}
    \subcaption{Inexact ($l=3$)}
  \end{minipage}
  \end{tabular}
  \caption{Designs after $3000$ iterations}
  \label{f_des_inexact}
\end{figure}

\begin{figure}[ht]
    \centering
    \includegraphics[width=8cm]{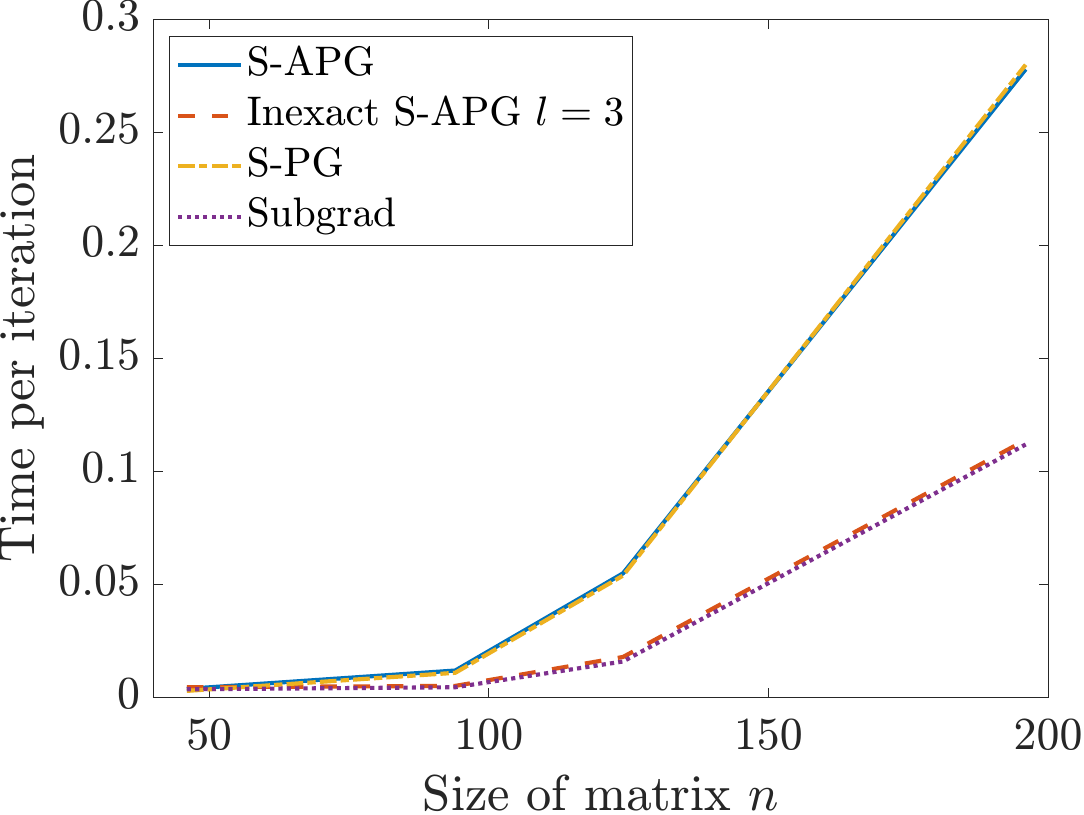}
    \caption{Computatinal costs per iteration for each dimension of the matrices $n$}
    \label{f_ite}
\end{figure}

To show the effectiveness of inexact smoothing, we compare S-APG and Inexact S-APG ($l=1,2,3$). Figure \ref{f_obj_inexact} shows the differences between objective values and the optimal value in 3000 iterations, Figure \ref{f_des_inexact} shows the truss designs after 3000 iterations, Figure \ref{f_ite} shows the computational cost per iteration for each size of the matrices $n$, and Table \ref{table1} shows the first to the third maximum eigenvalues\footnote{The maximum eigenvalue $\lambda_1\le 0$ and the minimal eigenfrequency $\omega_1\ge 0$ with the unit (1/s) has the following relation: $\lambda_1=-\omega_1^2$.} after $3000$ iterations including those of S-PG and Subgrad. We use MATLAB \texttt{eig} to compute all the eigenvalues in S-PG and S-APG, and MATLAB \texttt{eigs} (a sparse solver) to compute the largest $l$ eigenvalues in Inexact S-APG ($l=1,2,3$) and Subgrad ($l=1$). Note that the lines of S-APG and Inexact S-APG ($l=2,3$) overlap in Figure \ref{f_obj_inexact}. The lines of S-APG and S-PG and those of Inexact S-APG ($l=3$) and Subgrad also overlap in Figure \ref{f_ite}.

\begin{table}[ht]
\caption{Eigenvalues after $3000$ iterations}
\centering\small
\begin{tabular}{crrr}
\hline
Algorithm & $\lambda_1$ & $\lambda_2$ & $\lambda_3$ \\
\hline
S-APG & $-51.398$ & $-51.406$ & $-1094.580$ \\
Inexact S-APG ($l=1$) & $-51.300$ & $-51.406$ & $-7807.322$ \\
Inexact S-APG ($l=2$) & $-51.398$ & $-51.406$ & $-1097.965$ \\
Inexact S-APG ($l=3$) & $-51.398$ & $-51.406$ & $-1097.314$ \\
S-PG & $-50.561$ & $-50.961$ & $-24306.863$ \\
Subgrad & $-49.907$ & $-50.061$ & $-19417.493$ \\
\hline
\end{tabular}
\label{table1}
\end{table}

Figures \ref{f_des}(a), \ref{f_obj_inexact}, and \ref{f_des_inexact} show that Inexact S-APGs except for $l=1$ have the same performance as S-APG. Table \ref{table1} shows that the multiplicity of the maximum eigenvalue near the optimal solution is two, and thus Inexact S-APG with $l=1$ is not accurate enough. Figure \ref{f_ite} shows that Inexact S-APG can reduce the computational cost per iteration compared to S-APG, which computes all the eigenvalues.

From the above observation, we can expect that the inexact smoothing has enough accuracy when we set $l$ greater than the multiplicity of the maximum eigenvalues near the optimal solution. Although it is difficult to know the multiplicity in advance, it is much less than the size of the matrices $n$ in many cases.

\subsection{Comparison to the problem without the artificial lower bound}

\begin{figure}[ht]
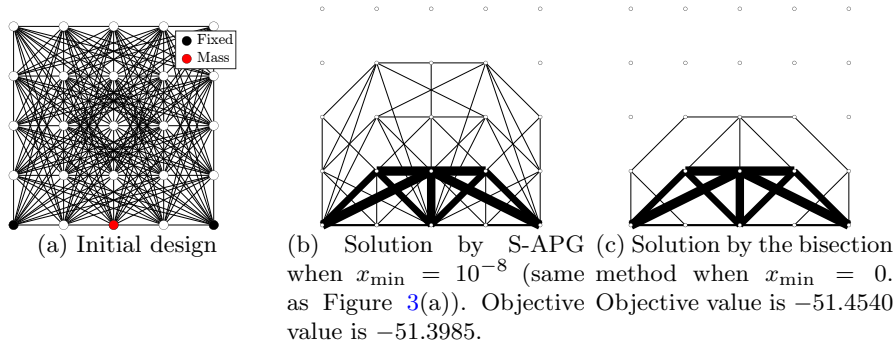

  \centering
  \begin{tabular}{ccc}
  \begin{minipage}[t]{0.3\hsize}
    \centering
    \includegraphics[width=3cm]{Fig_initial.pdf}
    \vspace{-2mm}
    \subcaption{Initial design}
  \end{minipage} &
  \hspace{-4mm}
  \begin{minipage}[t]{0.3\hsize}
    \centering
    \includegraphics[width=3cm]{Fig_sapg.pdf}
    \vspace{-2mm}
    \subcaption{Solution by S-APG when $x_\mathrm{min}=10^{-8}$ (same as Figure \ref{f_des}(a)). Objective value is $-51.3985$.}
  \end{minipage} &
  \hspace{-4mm}
  \begin{minipage}[t]{0.3\hsize}
    \centering
    \includegraphics[width=3cm]{Fig_bisection.pdf}
    \vspace{-2mm}
    \subcaption{Solution by the bisection method when $x_\mathrm{min}=0$. Objective value is $-51.4540$}
  \end{minipage}
  \end{tabular}
  \caption{Comparison of the solutions obtained by S-APG and the bisection method (case 1)}
  \label{f_bisection1}
\end{figure}

\begin{figure}[ht]
  \centering
  \begin{tabular}{ccc}
  \begin{minipage}[t]{0.3\hsize}
    \centering
    \includegraphics[width=3cm]{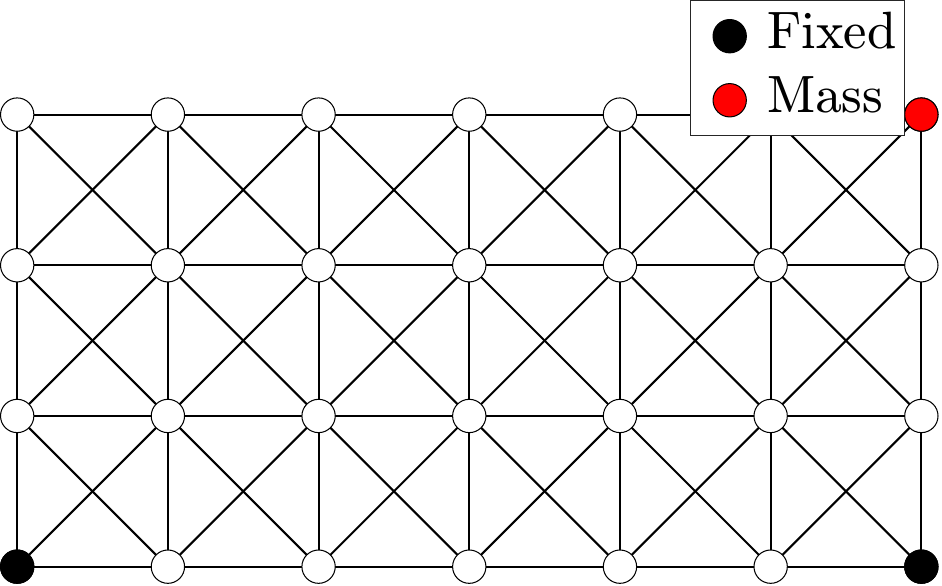}
    \vspace{-2mm}
    \subcaption{Initial design}
  \end{minipage} &
  \hspace{-4mm}
  \begin{minipage}[t]{0.3\hsize}
    \centering
    \includegraphics[width=3cm]{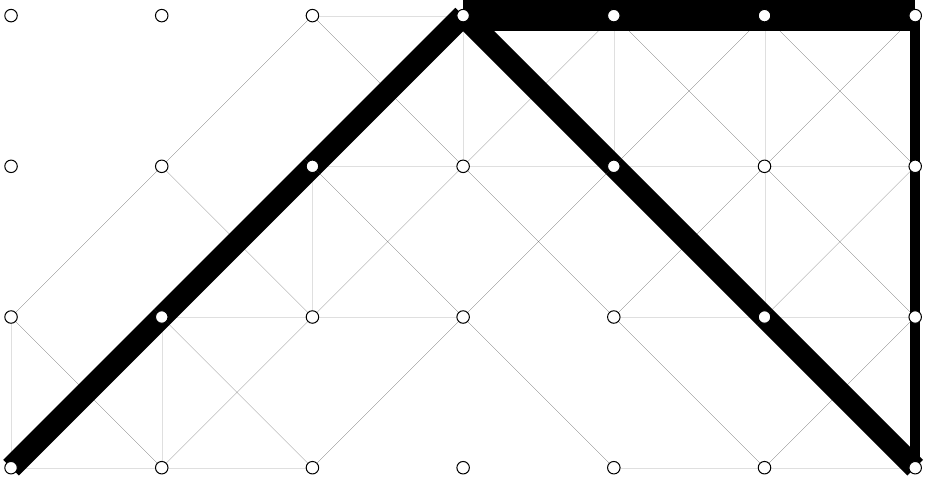}
    \vspace{-2mm}
    \subcaption{Solution by S-APG when $x_\mathrm{min}=10^{-8}$. Objective value is $-22.2205$}
  \end{minipage} &
  \hspace{-4mm}
  \begin{minipage}[t]{0.3\hsize}
    \centering
    \includegraphics[width=3cm]{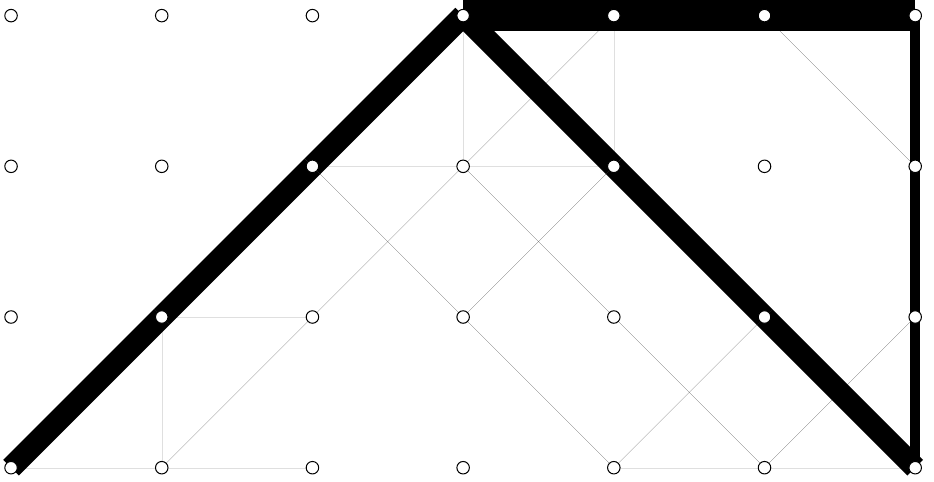}
    \vspace{-2mm}
    \subcaption{Solution by the bisection method when $x_\mathrm{min}=0$. Objective value is $-22.2216$}
  \end{minipage}
  \end{tabular}
  \caption{Comparison of the solutions obtained by S-APG and the bisection method (case 2)}
  \label{f_bisection2}
\end{figure}

We compare the solutions of the eigenfrequency optimization problem \eqref{maxeig} with the solutions of the problem \eqref{maxeig} without the artificial lower bound of the variables (cross-sectional areas of bars), namely the problem with $x_\mathrm{min}=0$. When $x_\mathrm{min}=0$, the maximum generalized eigenvalue can be discontinuous and our algorithms are not theoretically supported. However, it is still quasiconvex and the bisection method can be used to obtain the global optimal solution \cite{achtziger07siam}. See Appendices \ref{a_dis} and \ref{a_bis} for details. Note that the bisection method cannot be extended to problems when the feasibility subproblems are computationally costly, namely, when problems are large-scale or not quasiconvex (e.g., topology optimization of continua). Comparisons of the solutions obtained by S-APG and the bisection method in two different problem settings are shown in Figures \ref{f_bisection1} and \ref{f_bisection2}. The problem setting of Figure \ref{f_bisection2}(a) gives an example where the three maximum eigenvalues coincide near the optimal solution.


There are no significant differences between the solutions and their objective values obtained by S-APG and the bisection method shown in Figures \ref{f_bisection1} and \ref{f_bisection2} except for thin bars. This justifies the artificial lower bound of the variables to some extent.

\section{Conclusion}

In this paper, under some assumptions, we have investigated some properties of the maximum generalized eigenvalue: the Clarke subdifferential and pseudoconvexity. Moreover, algorithms to solve the maximum generalized eigenvalue minimization problem are considered. We have proved the convergence rate of the smoothing projected gradient method to the global optimum and proposed heuristic acceleration and inexact smoothing techniques to reduce practical computational costs.

Future work includes theoretical studies of acceleration, inexact smoothing, stepsize strategy, and stopping criteria, and extension to generalized eigenvalue optimization problems with more complicated constraints and symmetric-matrix-valued nonlinear functions. Also, a study of the discontinuous generalized eigenvalue with singular matrices (without the artificial lower bound of the optimization variables) is important. Although optimization problems involving generalized eigenvalues are an important class of problems in structural optimization, there are few theoretical studies. Further development of optimization theory and variational analysis of generalized eigenvalues is demanded.

\backmatter

\bmhead{Acknowledgments}

This research is part of the results of Value Exchange Engineering, a joint research project between R4D, Mercari, Inc. and RIISE. The work of the first author is partially supported by JSPS KAKENHI JP23KJ0383. The work of the third author is partially supported by JSPS KAKENHI JP19K15247. The work of the last author is partially supported by JSPS KAKENHI JP21K04351.

\begin{appendices}

\section{Generalized eigenvalue with singular matrices}
\label{a_dis}

In the eigenfrequency optimization \eqref{maxeig}, the matrices $K(\bm{x})$ and $M(\bm{x})+M_0$ can be singular when we set the lower bound of the variables as $x_{\mathrm{min}}=0$. When the matrices become singular (positive semidefinite) and $\bm{v}\in\mathrm{ker}K(\bm{x})\cap\mathrm{ker}M(\bm{x})+M_0$, $\lambda_i$ in the definition of generalized eigenvalue \eqref{ge} can take any values and is not well-defined. In \cite{achtziger07siam}, the extended definition of the minimum generalized eigenvalue of possibly singular matrices $(K(\bm{x}),M(\bm{x})+M_0)$ is introduced using the Rayleigh quotient:
\begin{align}
\lambda^{K,M}_n(\bm{x})\coloneq \underset{\bm{v}\notin\mathrm{ker}M(\bm{x})+M_0}{\inf}\frac{\bm{v}^\tp K(\bm{x})\bm{v}}{\bm{v}^\tp (M(\bm{x})+M_0)\bm{v}}.
\label{ge_semi}
\end{align}
We can define the maximum generalized eigenvalue as well $\lambda^{-K,M}_1(\bm{x})\coloneq-\lambda^{K,M}_n(\bm{x})$. Note that the properties of the stiffness and mass matrices, $K(\bm{x}),M(\bm{x})+M_0\succeq 0$ and $\mathrm{ker}(M(\bm{x})+M_0)\subseteq\mathrm{ker}K(\bm{x})\neq\rn$ for any $\bm{x}\in\rm_{\ge 0}\backslash\{\bm{0}\}$, are used in the definition and the proof of Proposition 2.3, and thus it is not directly extended to general symmetric matrices $A(\bm{x}),B(\bm{x})$. Proposition 2.3 in \cite{achtziger07siam} shows that the maximum generalized eigenvalue \eqref{ge_semi} is still quasiconvex\footnote{The proof of quasiconvexity in \cite{achtziger07siam} is incomplete. It relies on the fact that the maximum generalized eigenvalue \eqref{ge_semi} is written as a supremum of quasiconvex rational functions. However, whether a supremum of quasiconvex functions is quasiconvex is not clear if the index set depends on the variable $\bm{x}$ like \eqref{ge_semi} (consider, for example, a function $f_3:\mathbb{R}\to\mathbb{R}$ such that $f_3(x)=\max_{u\le 0}x^2+u$ if $x\le 1$ and $f_3(x)=\max_{u\le -1}x^2+u$ if $x>1$, which is not quasiconvex). Nevertheless, quasiconvexity can still be proved by the fact that the sublevel set is still written by $\{\bm{x}\in\rm \mid\lambda^{-K,M}_1(\bm{x})\le\alpha\}=\{\bm{x}\in\rm \mid K(\bm{x})+\alpha M(\bm{x}) \succeq 0\}$ (a direct consequence of Proposition 2.3(c) in \cite{achtziger07siam}), and it is convex for any $\alpha\in\mathbb{R}$.} on $\rm_{\ge 0}\backslash\{\bm{0}\}$, lower semicontinuous (continuous except on the boundary of $\rm_{\ge 0}\backslash\{\bm{0}\}$), and finite on $\rm_{\ge 0}\backslash\{\bm{0}\}$. Example 2.4 in \cite{achtziger07siam} gives an example where the minimum generalized eigenvalue is discontinuous on the boundary of $\rm_{\ge 0}\backslash\{\bm{0}\}$. 

By discontinuity, a global optimal solution of the minimization of $\lambda^{-K,M}_1(\bm{x})$ on $[x_{\mathrm{min}},\infty)^m,\ x_{\mathrm{min}}>0$, denoted by $\bm{x}^{*}$, is not necessarily close to a global optimal solution on $\rm_{\ge 0}\backslash\{\bm{0}\}$, denoted by $\bm{x}^{*}_0$. However, the quasiconvexity of the maximum generalized eigenvalue may restrict a possible position of $\bm{x}^{*}_0$. For example, $\bm{x}^{*}_0$ belongs to the sublevel set $\{\bm{x}\in\rm_{\ge 0}\backslash\{\bm{0}\}\mid\lambda^{-K,M}_1(\bm{x})\le\lambda_1^{-K,M}(\bm{x}^*)\}$ which is a convex set due to the quasiconvexity of the maximum generalized eigenvalue. Therefore, the line segment between $\bm{x}^{*}$ and $\bm{x}^{*}_0$ must belong to the sublevel set $\{\bm{x}\in\rm_{\ge 0}\backslash\{\bm{0}\}\mid\lambda^{-K,M}_1(\bm{x})\le\lambda_1^{-K,M}(\bm{x}^*)\}$, and a point in this line segment can belong to $[x_{\mathrm{min}},\infty)^m$ only if that point is also a global optimal solution on $[x_{\mathrm{min}},\infty)^m$. In particular, when $\bm{x}^{*}$ is the strict optimal solution on $[x_{\mathrm{min}},\infty)^m$, the line segment between $\bm{x}^{*}$ and $\bm{x}^{*}_0$ cannot intersect with $[x_{\mathrm{min}},\infty)^m\backslash\{\bm{x}^{*}\}$. This kind of property restrict a possible position of $\bm{x}^{*}_0$. Unfortunately, it is not easy to evaluate rigorously how close $\bm{x}^{*}$ and $\bm{x}^{*}_0$ are (or how similar their shapes are) and further theoretical studies are needed.

\section{A review of quasiconvex optimization algorithms}

In this section, we summarize implementations of existing algorithms for quasiconvex optimization in the maximum generalized eigenvalue minimization problem.

\subsection{Subgradient method}
\label{a_sub}

Since convex subdifferentials can be empty for quasiconvex functions, the quasiconvex subgradient method \cite{kiwiel01,konnov03,hu15} uses the closure of the Greenberg--Pierskalla (GP) subdifferential defined as follows.

\begin{definition}[Greenberg--Pierskalla (GP) subdifferential \cite{greenberg73}]
For a quasiconvex function $f:\rm\to\mathbb{R}$, the GP subdifferential at $\bm{x}\in\rm$ is defined as
\e{
\gpartial f(\bm{x})\coloneq\{\bm{g}\in\rm \mid \langle \bm{g},\bm{y}-\bm{x}\rangle <0,\ \forall\bm{y}\ \mathrm{s.t.}\ f(\bm{y}),<f(\bm{x})\},
\label{qcsub}
}
and each element of $\gpartial f(\bm{x})$ is called a GP subgradient. The closure of the GP subdifferential is often called the quasi-subdifferential (see \cite{hu15} for example).
\end{definition}

The GP subdifferential gives an optimality condition of a quasiconvex optimization problem. For any $\bm{x}\in\rm$, $\gpartial f(\bm{x})$ is nonempty, and $\bm{0}\in\gpartial f(\bm{x}^*)$, which is equivalent to $\gpartial f(\bm{x}^*)=\rn$, holds if and only if $\bm{x}^*\coloneq \mathrm{arg\,min}_x\ f(\bm{x})$.

The computation of a GP subgradient is impractical for some quasiconvex functions. However, for pseudoconvex functions, the definitions of pseudoconvexity and the GP subdifferential immediately lead to the fact that the Clarke subdifferential and the GP subdifferential have the inclusion $\cpartial f(\bm{x})\subseteq\gpartial f(\bm{x})$. The opposite inclusion is obviously false because $\gpartial f(\bm{x})$ is an unbounded convex cone. Note that the GP subdifferential can be defined for discontinuous functions, unlike the Clarke subdifferential.

Since the maximum generalized eigenvalue is pseudoconvex by Theorem \ref{t_pse}, the quasiconvex subgradient method \cite{kiwiel01,konnov03,hu15} for the maximum generalized eigenvalue minimization problem \eqref{p} becomes 
\begin{align}
& \bm{x}^{k+1}=\bm{x}^k-\alpha_k\frac{\bm{g}^k}{\|\bm{g}^k\|},\\
& \alpha_{k+1} = \alpha_0(k+1)^{-1/2},
\end{align}
where $\bm{g}^k\in\cpartial\lambda_1(A(\bm{x}^k),B(\bm{x}^k))$ is a Clarke subgradient.

\subsection{Bisection method}
\label{a_bis}

It is known that a global optimal solution of a quasiconvex optimization problem can be computed by the bisection method with a convex feasibility subproblem \cite{boyd04}. Consider a minimization problem of a quasiconvex function $f$. Set an estimate of the global optimal value $\lambda=(\underline{\lambda}+\overline{\lambda})/2\in\mathbb{R}$ where the interval $[\underline{\lambda},\overline{\lambda}]\subset\mathbb{R}$ is sufficiently large so that it contains the global optimal value, and solve the feasibility problem of the sublevel set $\{\bm{x}\in\rm\mid f(\bm{x})\le\lambda\}$, which is convex due to quasiconvexity of $f$. If it is feasible, the global optimal value is less than or equal to $\lambda$; otherwise, the global optimal value is greater than $\lambda$. Therefore, at each iteration of the bisection method, we can reduce the size of the interval $[\underline{\lambda},\overline{\lambda}]$, containing the global optimal value, to half.

The bisection method for the maximum generalized eigenvalue minimization problem \cite{achtziger07siam} solves the feasibility problem
\e{
\text{Find }\bm{x}\in S\text{ s.t. }A(\bm{x})-\lambda B(\bm{x})\preceq 0
\label{fea1}
}
for fixed $\lambda\in\mathbb{R}$, which can be solved by a standard linear semidefinite programming solver (we use SDPT3 of CVX in our numerical experiments). However, the feasibility problem becomes hard to solve when $\lambda$ is very close to the global optimal value because the feasible set becomes very small. Therefore, we propose a modification; instead of \eqref{fea1}, we solve the minimization problem
\e{
\begin{aligned}
& \underset{\bm{x}\in S,\ z\in\mathbb{R}}{\mathrm{Minimize}} & & z \\
& \mathrm{subject\ to} & & A(\bm{x})-\lambda B(\bm{x})-zI\preceq 0,
\label{fea2}
\end{aligned}
}
with an auxiliary variable $z\in\mathbb{R}$. The problem \eqref{fea2} is always feasible, and the condition that the optimal value of \eqref{fea2} is nonpositive is equivalent to the feasibility of \eqref{fea1}.




\end{appendices}


\bibliography{ref}

\end{document}